\documentclass[a4paper]{amsart}

\usepackage[utf8]{inputenc}
\usepackage{amssymb,amsfonts,amsmath}
\usepackage{enumerate}
\usepackage{color}
\usepackage{mathrsfs}
\usepackage{graphicx}
\usepackage{verbatim}

\newtheorem{theorem}{Theorem}[section]

\newtheorem{corollary}[theorem]{Corollary}
\newtheorem{lemma}[theorem]{Lemma}
\newtheorem{proposition}[theorem]{Proposition}
\newtheorem{question}[theorem]{Question}

\theoremstyle{definition}
\newtheorem{definition}[theorem]{Definition}

\newtheorem{notation}[theorem]{Notation}

\theoremstyle{remark}
\newtheorem{remark}[theorem]{Remark}

\newtheorem{example}[theorem]{Example}

\newcommand{\N}{\mathbb{N}}

\newcommand{\R}{\mathbb{R}}
\newcommand{\Z}{\mathbb{Z}}

\DeclareMathOperator{\CAT}{CAT}
\DeclareMathOperator{\Hd}{Hd}
\DeclareMathOperator{\MCG}{MCG}
\DeclareMathOperator{\id}{id}
\DeclareMathOperator{\im}{im}
\DeclareMathOperator{\supp}{supp}
\DeclareMathOperator{\Out}{Out}
\newcommand{\abs}[1]{\vert #1 \vert}

\newcommand{\defeq}{\mathrel{\mathop{:}}=}

\begin{document}

\title[Relative growth rates in acylindrically hyperbolic groups]{The relative exponential growth rate of subgroups of acylindrically hyperbolic groups}


\date{\today}		   	
		   			


\author[E.~Schesler]{Eduard Schesler}
\address{Fakult\"{a}t f\"{u}r Mathematik und Informatik, FernUniversit\"{a}t in Hagen, 58084 Hagen, Germany}
\email{eduard.schesler@fernuni-hagen.de}

\thanks{The author was partially supported by the DFG grant WI 4079/4
within the SPP 2026 Geometry at infinity.}

\begin{abstract}
We introduce a new invariant of finitely generated groups, the ambiguity function, and prove that every finitely generated acylindrically hyperbolic group has a linearly bounded ambiguity function.
We use this result to prove that the relative exponential growth rate $\lim \limits_{n \rightarrow \infty} \sqrt[n]{\vert B^{X}_H(n) \vert}$ of a subgroup $H$ of an acylindrically hyperbolic group
$G$ exists with respect to every finite generating set $X$ of $G$, if $H$ contains a loxodromic element of $G$.
Further we prove that the relative exponential growth rate of every finitely generated subgroup $H$ of a right-angled Artin group $A_{\Gamma}$ exists with respect to every finite generating set of $A_{\Gamma}$.
\end{abstract}

\maketitle

\section{Introduction}

Given a finitely generated group $G$ and a finite generating set $X$ of $G$, the \emph{growth function} $\beta^X_G \colon \N \rightarrow \N$ of $G$ with respect to $X$ is defined by $\beta^X_G(n) = \abs{B^X_G(n)}$, where $B^X_G(n)$ denotes the set of all elements of $G$ that are represented by words of length $\leq n$ in the generators of $X$ and $X^{-1}$.
By definition, $G$ has \emph{exponential growth} if there is a constant $c>1$ such that $\beta^X_G(n) > c^n$ for every $n \in \N$.
Using the submultiplicative property
\[
\beta^X_G(m+n) \leq \beta^X_G(m)\beta^X_G(n),
\]
Milnor proved in~\cite{Milnor68} that the limit $a \defeq \lim \limits_{n \rightarrow \infty} \sqrt[n]{\beta^X_G(n)}$ exists for every finitely generated group $G$ and every finite generating set $X$ of $G$.
Thus for every $\varepsilon > 0$ there is a number $N \in \N$ such that for every $n \geq N$ we have
\[
(a-\varepsilon)^n \leq \beta^X_G(n) \leq (a+\varepsilon)^n.
\]
If the limit $a$ is larger than $1$, this justifies the terminology of $G$ having exponential growth and that $a$ is called the \emph{exponential growth rate}, or just \emph{growth rate}, of $G$ with respect to $X$.
In the following
In this paper we study a relative analogue of the growth function.
For every (not necessarily finitely generated) subgroup $H \leq G$ and every $n \in \N$ we consider the subset
\[
B^X_H(n) \defeq \{h \in H \colon \abs{h}_X \leq n \} = B^X_G(n) \cap H
\]
of $B^X_G(n)$.
The \emph{relative growth function} $\beta^X_H \colon \N \rightarrow \N$ of $H$ with respect to $X$ is defined by $\beta^X_H(n) = \abs{B^X_H(n)}$.

Again we can ask for the existence of the limit $\lim \limits_{n \rightarrow \infty} \sqrt[n]{\beta^X_H(n)}$.
In the case where this limit exists we will say that the relative (exponential) growth rate $\lim \limits_{n \rightarrow \infty} \sqrt[n]{\beta^X_H(n)}$ of $H$ exists with respect to the generating set $X$ of $G$.
The first results in this direction were obtained by Cohen~\cite{Cohen82} and Grigorchuk~\cite{Grigorchuk80b}.
Given a finitely generated free group $F$ with basis $X$ and a normal subgroup $N \unlhd F$, Cohen and Grigorchuk independently proved that the relative growth rate of $N$ exists with respect to $X$ and that the quotient $F/N$ is amenable if and only if the growth rates $\lim \limits_{n \rightarrow \infty} \sqrt[n]{\beta^X_F(n)}$ and $\lim \limits_{n \rightarrow \infty} \sqrt[n]{\beta^X_N(n)}$ coincide.
Note that in this case
\[
\lim \limits_{n \rightarrow \infty} \sqrt[n]{\beta^X_N(n)} = 2 \abs{X} -1.
\]
Later it turned out that the assumption of being normal is not necessary for the relative growth rate to exist.
In~\cite{Olshanskii17} Olshanskii showed that the relative growth rate of every subgroup of a finitely generated free group $F$ exists with respect to a free generating set of $F$.
Furthermore Olshanskii pointed out that the relative growth rate does not exist in general.
More precisely he showed that there is a finitely generated solvable group $G$ with a finite generating set $X$ and an infinite cyclic subgroup $H \leq G$ such that
\[
\liminf \limits_{n \rightarrow \infty} \frac{\beta^X_H(n)}{n^{1+\varepsilon}} = 0
\hspace*{3mm}\forall\varepsilon > 0 \hspace*{3mm} \text{and}\hspace*{3mm}
\limsup \limits_{n \rightarrow \infty} \sqrt[n]{\beta^X_H(n)} > 1.
\]
Meanwhile it turned out that many results about relative growth rates of subgroups in free groups have analogues for hyperbolic groups.
Coulon, Dal'Bo and Sambusetti generalized the amenability criterion of Cohen and Grigorchuk to quotients of hyperbolic groups~\cite{CoulonDalBoSambusetti18}.
For a hyperbolic group $G$ with finite generating set $X$ and a normal subgroup $N \unlhd G$ they showed that the quotient $G/N$ is amenable if and only if
\[
\limsup \limits_{n \rightarrow \infty} \sqrt[n]{\beta^X_G(n)} = 
\limsup \limits_{n \rightarrow \infty} \sqrt[n]{\beta^X_N(n)}.
\]
Further it was proven by Dahmani, Futer, and Wise that the relative growth rate of a quasi-convex subgroup $H$ of a hyperbolic group $G$ is strictly less than the relative growth rate of $G$ (see~\cite{DahmaniFuterWise19}).
More recently, this result was extended by Cordes, Russell, Spriano, and Zalloum to more general classes of groups (see~\cite{CordesRussellSprianoZalloum20}).
It is now natural to ask whether the existence result of Olshanskii can be generalized to hyperbolic groups, which will allow us to replace $\limsup \limits_{n \rightarrow \infty}$ by $\lim \limits_{n \rightarrow \infty}$ when speaking about relative growth rates.

\begin{question}\label{QuestionExistenceForHyperbolic}
Let $G$ be a hyperbolic group, let $X$ be a finite generating set of $G$ and let $H$ be an arbitrary subgroup of $G$.
Does the relative growth rate $\lim \limits_{n \rightarrow \infty} \sqrt[n]{\beta^X_H(n)}$ exist?
\end{question}

This might seem unlikely since we know from the Rips construction~\cite{Rips82} that there are hyperbolic groups with extremely distorted and bad-behaved subgroups.
Further it was pointed out by Sharp in~\cite{Sharp98} that there are hyperbolic groups $G$ with finitely generated normal subgroups $N$ such that the relative growth function $\beta^{X}_N$ can be very wild for some (and hence any) finite generating set $X$ of $G$.
By analyzing the asymptotic behaviour of $\beta^{X}_N$, Sharp showed that the formal power series $\sum \limits_{n = 0}^{\infty} \beta^{X}_N(n) z^n$ is not rational in general.
This is in contrast to the case of finitely generated free groups $F(X)$ for which Grigorchuk~\cite{Grigorchuk80b} showed that $\sum \limits_{n = 0}^{\infty} \beta^{X}_H(n) z^n$ is rational for every finitely generated subgroup $H \leq F(X)$.
Nevertheless we will see that, even in the larger context of acylindrically hyperbolic groups, the relative growth rate of subgroups exists under mild conditions.

\begin{theorem}\label{IntroMainTheorem}
Let $G$ be a finitely generated acylindrically hyperbolic group and let $H$ be a subgroup of $G$.
If $H$ contains a generalized loxodromic element of $G$, then the relative growth rate of $H$ exists with respect to every finite generating set of $G$.
\end{theorem}



From this it will be easy to derive an affirmative answer to Question~\ref{QuestionExistenceForHyperbolic}.

\begin{corollary}\label{IntroCorollaryExistForHypGrps}
Let $G$ be a hyperbolic group and let $H \leq G$ be an arbitrary subgroup.
The relative growth rate of $H$ exists with respect to every finite generating set of $G$.
\end{corollary}


In order to prove Theorem~\ref{IntroMainTheorem}, we introduce the (relative) ambiguity function (see Definition~\ref{DefinitionRelativeAmbiguity}), which is a new invariant of finitely generated groups.
Coarsely speaking, the ambiguity function of a group $G$ counts the number of ways an element $g \in G$ can be written as a product $g = h x_{h,k} k$ of elements $h, x_{h,k}, k \in G$ of prescribed lengths, where the length of $x_{h,k}$ is uniformly bounded.
We will prove the following.

\begin{theorem}\label{IntroTheoremAmbuguityProducts}
Let $G = \prod \limits_{i = 1}^{n} G_i$ be a product of finitely generated acylindrically hyperbolic groups $G_i$.
Then $G$ has an ambiguity function that is bounded above by a polynomial of degree $n$.
\end{theorem}

As an application we will study relative growth rates in right-angled Artin groups, which are known to decompose as direct products of acylindrically hyperbolic groups and infinite cyclic groups.

\begin{theorem}\label{IntroTheoremExistenceForRAAGs}
Let $A_{\Gamma}$ be a right-angled Artin group and let $H \leq A_{\Gamma}$ be a finitely generated subgroup.
Then the relative growth rate of $H$ exists with respect to every finite generating set of $A_{\Gamma}$.
\end{theorem}

\smallskip

\subsection*{The paper is organized as follows.}
In Section~\ref{SectionCalculus} we will introduce the (relative) ambiguity function of a subgroup of a finitely generated group and we will see how this function relates to the existence of relative growth rates.
In Section~\ref{SectionHyperbolicSpaces} we recall some basic concepts about $\delta$-hyperbolic spaces and prove some auxiliary results.
The necessary terminology and some well-known results about acylindrically hyperbolic groups will be introduced in Section~\ref{SectionAcylindricallyHypGrps}.
In Section~\ref{SectionExistence} we will study relative ambiguity functions in acylindrically hyperbolic groups and we will prove Theorem~\ref{IntroMainTheorem}.
From this we will derive a positive answer to Question~\ref{QuestionExistenceForHyperbolic}.
In Section~\ref{SectionProducts} we will extend the study of relative ambiguity functions to products of acylindrically hyperbolic groups.
We will prove Theorem~\ref{IntroTheoremAmbuguityProducts}, which will allow us to show that an analogue of Theorem~\ref{IntroMainTheorem} holds for products of acylindrically hyperbolic groups.
As an application, we will be able to study relative growth rates in right-angled Artin groups and to prove Theorem~\ref{IntroTheoremExistenceForRAAGs} in Section~\ref{Applications}.

\section{Weak Supermultiplicativity}\label{SectionCalculus}

The aim of this section is to formulate sufficient conditions for the relative growth rate to exist.
In order to do so, we will introduce the relative ambiguity function and combine it with the well-known distortion function.

\begin{definition}
A function $f \colon \N \rightarrow \N$ is called \emph{sublinear} if $\lim \limits_{n \rightarrow \infty} \frac{f(n)}{n} = 0$ and \emph{subexponential} if $\lim \limits_{n \rightarrow \infty} \frac{f(n)}{a^n} = 0$ for every $a > 1$.
\end{definition}

\begin{definition}\label{DefinitionDistortion} 
Let $G$ be a finitely generated group with a finite generating set~$X$.
Let $H$ be a finitely generated subgroup of $G$ with a finite generating set $Y$.
The \emph{distortion function} of $H$ in $G$ with respect to $X$ and $Y$ is defined by
\[
\Delta_H^G \colon \N \rightarrow \N,\
n \mapsto \max\{\abs{h}_Y \colon h \in H, \ \abs{h}_X \leq n \}.
\]
\end{definition}

We will be especially interested in the cases where $\Delta_H^G$ is bounded above by a linear, polynomial or subexponential function.
Note that these properties do not depend on the choices of the finite generating sets of $H$ and $G$.
\\\\Suppose now that $H$ is an arbitrary (not necessarily finitely generated) subgroup of a finitely generated group $G$ and let $X$ be a finite generating set of $G$.
The following maps, indexed by $m,n\in \N$, are given by restricting the group multiplication in $H$ to relative balls of radius $m$ and $n \colon$
\[
\mu_{m,n} \colon B_H^X(m)\times B_H^X(n)\rightarrow B_H^X(m+n), \hspace*{2mm} (h,k) \mapsto hk.
\]
The main difficulty in showing that the relative growth rate of $H$ exists with respect to $X$ is that $\mu_{m,n}$ may fail to be surjective.
In this case we lose the submultiplicative property
\[
\beta^X_H(m+n) \leq \beta^X_H(m)\beta^X_H(n),
\]
which was used in~\cite{Milnor68} to prove that $\lim \limits_{n \rightarrow \infty} \sqrt[n]{\beta^X_H(n)}$ exists in the case where $H=G$.
Our goal in the following will be to replace the submultiplicativity by a weak form of supermultiplicativity.
More precisely, we want to establish inequalities of the form
\[
\beta^X_H(m) \beta^X_H(n) \leq \varepsilon(n) \beta^X_H(m+n+c)
\]
for some constant $c$ and some subexponential function $\varepsilon \colon \N \rightarrow \N$ that does not depend on the radius $m$.
A naive approach to obtain such an inequality is given by setting $c = 0$ and to define
\[
\varepsilon(n) \defeq \sup \limits_{m \in \N} \ \max \limits_{g \in B_H^X(m+n)} \abs{\mu_{m,n}^{-1}(g)}
\]
as the maximal cardinality of the fibers $\mu_{m,n}^{-1}(g)$.
Note that $\varepsilon$ is well-defined since $\abs{\mu_{m,n}^{-1}(g)} \leq \beta^X_H(n)$ for every $n \in \N$.
However, if $G$ is any finitely generated group of exponential growth with a finite generating set $X$ and $H=G$, we see that
\[
\abs{\mu_{n,n}^{-1}(1)} = \abs{\{(g,g^{-1}) \in B_{G}^X(n) \times B_{G}^X(n) : g \in B_{G}^X(n) \}}
= \beta^X_G(n)
\]
growth exponenitally with respect to $n$.
In~\cite[Theorem~1.5]{Olshanskii17} Olshanski was able to overcome this problem in the case where $G = F(X)$ is a free group with basis $X$ and $H$ is a subgroup of $F(X)$.
In order to reduce the number of elements in the fibers $\mu_{m,n}^{-1}(g)$, he modified the maps $\mu_{m,n}$ by inserting an appropriate element $x_{h,k} \in H$ between every pair of elements $h,k \in H$ before multiplying them.
The following definition is an attempt to capture this idea.

\begin{definition}\label{DefinitionRelativeAmbiguity}
Let $G$ be a finitely generated group with a finite generating set $X$.
Let $H \leq G$ be an arbitrary subgroup.
A function $\Phi \colon H \times H \rightarrow H$ is a \emph{concatenation map} of $H$ in $G$ if there is a constant $\varepsilon \in \N$ such that
\[
\Phi(h,k) = h x_{h,k} k \text{ with } x_{h,k} \in B^X_H(\varepsilon) \text{ for all } h,k \in H.
\]
We will refer to $B^X_H(\varepsilon)$ as a set of \emph{connecting pieces} for $\Phi$.
For all $m,n \in \N$ let $\Phi_{m,n}$ denote the restriction of $\Phi$ to $B^X_H(m) \times B^X_H(n)$.
A function $f \colon \N \rightarrow \N$ is a \emph{relative ambiguity function} for $\Phi$ if
\[
\abs{\Phi_{m,n}^{-1}(g)} \leq f(n) \text{ for all } m,n \in \N \text{ and } g \in H.
\]
In this case we also say that $f$ is a \emph{relative ambiguity function} of $H$ in $G$.
If $H = G$ we just say that $f$ is an \emph{ambiguity function} of $G$.

\end{definition}


The first thing to observe is that, up to precomposing with a linear function, the relative ambiguity function does not depend on the choice of the finite generating set of $G$.

\begin{lemma}\label{LemmaIndepOfGenSet}
Let $G$ be a finitely generated group and let $H$ be a subgroup of $G$.
Let $X$ and $Y$ be finite generating sets of $G$.
Then there is a constant $c > 0$ such that for every relative ambiguity function $f$ of $H$ in $G$ with respect to $X$ the function
$n \mapsto f(cn)$ is a relative ambiguity function of $H$ in $G$ with respect to $Y$.
\end{lemma}
\begin{proof}
Suppose that $f$ corresponds to a concatenation map $\Phi \colon H \times H \rightarrow H$ with $B^X_H(\varepsilon)$ as a set of connecting pieces.
By setting $c \defeq \max \limits_{y \in Y} \abs{y}_{X}$, we obtain the inclusions
\[
\iota_{m,n} \colon B^Y_H(m) \times B^Y_H(n) \subseteq B^X_H(cm) \times B^X_H(cn) \text{ for all } m,n \in \N.
\]
Note that the restriction of $\Phi$ to $B^Y_H(m) \times B^Y_H(n)$ can be written as $\Phi_{m,n}' \defeq \Phi_{cm,cn} \circ \iota_{m,n}$
and that the definition of $f$ gives us
\[
\abs{\Phi_{m,n}'^{-1}(w)} = \abs{\iota_{m,n}^{-1}(\Phi_{cm,cn}^{-1}(w))}
\leq \abs{\Phi_{cm,cn}^{-1}(w)} \leq f(cn)
\]
for every $w \in H$.
Thus we see that $n \mapsto f(cn)$ defines a relative ambiguity function of $H$ in $G$ with respect to $Y$.
\end{proof}

Note that Lemma~\ref{LemmaIndepOfGenSet} implies in particular that the properties of $H$ to have a relative ambiguity function in $G$ that is bounded below (resp. above) by a linear, polynomial, or subexponential function does not depend on the choice of the finite generating set of $G$.

\medskip

The following example is a special case of the construction in the proof of~\cite[Theorem~1.5]{Olshanskii17}.

\begin{example}\label{ExampleAmbiguityFreeGroup}
We consider a finitely generated non abelian free group $F = F(X)$ with basis $X$.
Note that for each pair of reduced words $(u,v) \in F \times F$ there is a letter $x_{u,v} \in X \cup X^{-1}$ such that $u x_{u,v} v$ is reduced.
This gives us a concatenation map
\[
\Phi \colon F \times F \rightarrow F,\ (u,v) \mapsto u x_{u,v} v
\]
with $B^X_F(1)$ as a set of connecting pieces.
In order to define an ambiguity function for $F$, we need to estimate the cardinality of each fiber $\Phi_{m,n}^{-1}(w) \subset B^X_{F}(m) \times B^X_{F}(n)$ where $w \in F$.
Thus let $(u,v) \in B^X_{F}(m) \times B^X_{F}(n)$ be a pair of (reduced) elements such that $ux_{u,v}v = w$.
Since $ux_{u,v}v$ is reduced, it follows that $u$ and $v$ are completely determined by their length and we obtain $\abs{\Phi_{m,n}^{-1}(w)} \leq \min(m,n)$.
In particular we see that the identity $\id \colon \N \rightarrow \N$ is an ambiguity function of $F$.
\end{example}

In Section~\ref{SectionExistence} we will see that the existence of a linear ambiguity function is a typical consequence in the presence of negative curvature.

\begin{example}\label{ExampleAmbiguityAbelianGroup}
Consider the free abelian group $\Z^2$ of rank $2$ and its standard generating set $X = \{(1,0),(0,1)\}$.
We claim that $\Z^2$ does not posses a linearly bounded ambiguity function.
To see this let
\[
\Phi \colon \Z^2 \times \Z^2 \rightarrow \Z^2,\ (u,v) \mapsto u + g_{u,v} + v
\]
be a concatenation map with $B^X_{\Z^2}(\varepsilon)$ as a set of connecting pieces for some $\varepsilon \in \N$. 
Note that we have
\[
\Phi(u,-u) = u + g_{u,-u} - u = g_{u,-u} \in B_{\Z^2}^X(\varepsilon) \text{ for every } u \in \Z^2.
\]
In particular we see that for every $n \in \N$ there is an element $w \in B^X_{\Z^2}(\varepsilon)$ such that the number of elements in the fiber of $w$ with respect to $\Phi_{n,n}$ satisfies
\[
\abs{\Phi_{n,n}^{-1}(w)} \geq \frac{\beta_{\Z^2}^X(n)}{\beta_{\Z^2}^X(\varepsilon)} \geq \frac{n^2}{\beta_{\Z^2}^X(\varepsilon)}.
\]
Hence every ambiguity function of $\Z^2$ grows quadratically.
\end{example}


%




\begin{lemma}\label{LemmaSubexpUndistorted}
Let $G$ be a finitely generated group with a finite generating set $X$ and let $H \leq G$ be a finitely generated, subexponentially distorted subgroup.
If $H$ possesses an ambiguity function that is bounded above by a polynomial,
then there is a subexponential function $\delta \colon \N \rightarrow \N$ and a constant $\varepsilon \in \N$ such that
\[
\beta^{X}_H(m) \beta^{X}_H(n) \leq \delta(n) \cdot \beta^{X}_H(m+n+\varepsilon)
\hspace*{2mm}\text{ for all } m,n \in \N.
\]
\end{lemma}
\begin{proof}
Let $Y$ be a finite generating set of $H$ and let $f$ be a subexponential function bounding the distortion of $H$ in $G$ with respect to $X$ and $Y$.
This gives us a well defined inclusions
\[
\iota_{m,n} \colon B^{X}_H(m) \times B^{X}_H(n) \rightarrow B^{Y}_H(f(m)) \times B^{Y}_H(f(n)).
\]
Without loss of generality we may assume that $Y$ is a subset of $X$.
Let $P$ be a polynomially bounded ambiguity function of $H$ with respect to $Y$, some concatenation map $\Phi$, and a set $B^{Y}_H(\varepsilon)$ of connecting pieces.
Note that $\Phi$ can be restricted to the maps
\[
\Phi_{f(m),f(n)} \colon B^{Y}_H(f(m)) \times B^{Y}_H(f(n)) \rightarrow B^{Y}_H(f(m) + f(n) + \varepsilon),
\ (h,k) \mapsto h x_{h,k} k
\]
for all $m,n \in \N$, where $x_{h,k} \in B^{Y}_H(\varepsilon)$ and $\abs{\Phi_{f(m),f(n)}^{-1}(z)} \leq P(f(n))$ for every $z \in B^{Y}_H(f(m) + f(n) + \varepsilon)$.
By precomposing these restrictions with $\iota_{m,n}$, we obtain maps
\[
B^{X}_H(m) \times B^{X}_H(n) \rightarrow B^{Y}_H(f(m) + f(n) + \varepsilon) \text{ for all } m,n \in \N
\]
whose fibers are bounded by $P(f(n))$.
In particular we get
\[
\beta^{X}_H(m) \beta^{X}_H(n) \leq P(f(n)) \cdot \beta^{Y}_H(m+n+\varepsilon) \text{ for all } m,n \in \N.
\]
Since $Y \subset X$ implies $\beta^{Y}_H(n) \leq \beta^{X}_H(n)$ for every $n \in \N$, the lemma follows from the easy observation that $\delta \defeq P \circ f$ is a subexponential function.
\end{proof}



In the case where $H$ is \emph{undistorted} in $G$, i.e.\ where $\Delta_H^G$ is bounded above by a linear function, we get the following analogue of Lemma~\ref{LemmaSubexpUndistorted}.


\begin{lemma}\label{LemmaNewSubexponential}
Let $G$ be a finitely generated group and let $H \leq G$ be a finitely generated, undistorted subgroup of $G$.
If $H$ possesses a subexponential ambiguity function, then there is a subexponential function $\delta \colon \N \rightarrow \N$ and a constant $\varepsilon \geq 0$ such that
\[
\beta^{X}_H(m) \beta^{X}_H(n) \leq \delta(n) \cdot \beta^{X}_H(m+n+\varepsilon) \text{ for all } m,n \in \N.
\]
\end{lemma}
\begin{proof}
The proof is essentially the same as the proof of Lemma~\ref{LemmaSubexpUndistorted}.
One only has to replace the fact that the composition $P \circ f$ of a polynomial $P$ and a subexponential function $f$ is subexponential by the fact that $f \circ L$ is subexponential if $L$ is a linear function.
\end{proof}

The following lemma provides us with an alternative definition of a subexponential function, which will be easier to apply for our purposes.

\begin{lemma}\label{LemmaNewSubexpPhrasing}
A function $f \colon \N \rightarrow \N$ is subexponential if and only if $\lim \limits_{n \rightarrow \infty} f(n)^{\frac{1}{n}} = 1$.
\end{lemma}
\begin{proof}
First we assume that $f$ is subexponential.
Then for every $a > 1$ and every $\varepsilon > 0$ we have $\frac{f(n)}{a^n} \leq \varepsilon$ for almost every $n \in \N$.
By choosing $\varepsilon \leq 1$, we obtain $f(n)^{\frac{1}{n}} \leq a$ for almost every $n \in \N$.
Since we can choose $a > 1$ arbitrarily close to $1$, we get $\lim \limits_{n \rightarrow \infty} f(n)^{\frac{1}{n}} = 1$.
Suppose now that $\lim \limits_{n \rightarrow \infty} f(n)^{\frac{1}{n}} = 1$.
Thus, for every $\varepsilon > 0$ we have $f(n) \leq (1 + \varepsilon)^n$ for almost every $n \in \N$.
Let $\delta > 0$ and let $a \defeq 1 + \delta$.
By setting $\varepsilon = \frac{\delta}{2}$, we obtain
\[
\limsup \limits_{n \rightarrow \infty} \frac{f(n)}{a^n}
\leq \limsup \limits_{n \rightarrow \infty} \frac{(1 + \varepsilon)^n}{a^n}
= \limsup \limits_{n \rightarrow \infty} \Big( \frac{1 + \frac{\delta}{2}}{1 + \delta} \Big)^n
= 0
\]
and therefore $\lim \limits_{n \rightarrow \infty} \frac{f(n)}{a^n} = 0$.
\end{proof}

The crucial analytical ingredient for the proof of the existence of some relative growth rates is provided by the following proposition, which is an adaptation of the analytical part of the proof of Lemma $1.6$ in~\cite{Olshanskii17} for our purposes.
For convenience we include a proof here.

\begin{proposition}\label{PropositionCalculusOptimized}
Let $B,C > 0$ be two real numbers, $\delta \colon \N \rightarrow \N$ a subexponential function, and $l \colon \N \rightarrow \N$ a sublinear function.
If $f \colon \N \rightarrow \N$ is a monotonically increasing function with the properties
\begin{enumerate}
\item[\rm{1)}] $f(m)f(n)\leq \delta(n)f(m+n+l(n))$ for all $n \geq C$
and all $m \in \N$ and \vspace{2mm}
\item[\rm{2)}] $1 \leq f(n) \leq B^n$ for all $n \in \N$,
\end{enumerate}
then the limit $\lim \limits_{n \rightarrow \infty}{f(n)^{\frac{1}{n}}}$ exists.
\end{proposition}
\begin{proof}
First we observe that $\limsup \limits_{n\rightarrow \infty}{f(n)^{\frac{1}{n}}}$ exists due to property 2).
Without loss of generality we may assume that $a \defeq \limsup \limits_{n\rightarrow \infty}{f(n)^{\frac{1}{n}}} > 1$.
Indeed, since $f(n) \geq 1$ for every $n \in \N$, we would otherwise obtain
\[
1 \leq \liminf \limits_{n\rightarrow \infty}{f(n)^{\frac{1}{n}}}
\leq \limsup \limits_{n\rightarrow \infty}{f(n)^{\frac{1}{n}}} \leq 1
\]
in which case $\lim \limits_{n\rightarrow \infty}{f(n)^{\frac{1}{n}}}$ clearly exists.
Let $\varepsilon \in (0,a-1)$ and let $a_n \defeq f(n)^{\frac{1}{n}}$.
By Lemma~\ref{LemmaNewSubexpPhrasing} we have $\lim \limits_{n \rightarrow \infty}{\delta(n)^{\frac{1}{n}}} = 1$ and therefore $\lim \limits_{n \rightarrow \infty}{\delta(n)^{\frac{-1}{n}}} = 1$.
Since $l$ is sublinear, we further have $\lim \limits_{n \rightarrow \infty}{\frac{n+l(n)}{n}} = 1$.
Thus, there is a natural number $N \geq C$ such that
\begin{equation}\label{EquationEpsOverThree}
|a_{N} \delta(N)^{\frac{-1}{N}}-a| < \frac{\varepsilon}{3} \text{ and }
\Big(a-\frac{2\varepsilon}{3}\Big)^{\frac{N+l(N)}{N}} < a-\frac{\varepsilon}{3}.
\end{equation}
In the following, a natural number $n$ will be written as
$n = q(N+l(N))+r$ with $r < N+l(N)$ and $q \in \N$.
Since $f$ is monotonically increasing and
\begin{equation}\label{EquationNandlN}
\Big(a-\frac{2\varepsilon}{3}\Big)^{n-N-l(N)} > (a-\varepsilon)^{n}
\end{equation}
holds for almost every $n \in \N$, an iterative application of property 1) gives us
\begin{center}
\begin{tabular}{lll}
$f(n)$ & $\geq f(q(N+l(N)))$ & $\geq \delta(N)^{-1}f(N)f((q-1)(N+l(N)))$\\[2ex]
       & $ \geq \hspace{1cm} \dots $   & $\geq \delta(N)^{-q+1}f(N)^{q-1} f(N+l(N))$\\[2ex]
       & $\geq \delta(N)^{-q}f(N)^q$. &
\end{tabular}
\end{center}
Since $a_n = f(n)^{\frac{1}{n}}$ by definition, we therefore obtain
\begin{equation}\label{EquationSubmult}
f(n) \geq \delta(N)^{-q}f(N)^q = (\delta(N)^{\frac{-1}{N}} a_N )^{Nq}.
\end{equation}
By combining the above inequalities, we get
\begin{center}
\begin{tabular}{llll}
$f(n)$ & $\stackrel{\eqref{EquationSubmult}}{\geq} \Big(a_N \delta(N)^{\frac{-1}{N}}\Big)^{Nq}$ &$ \stackrel{\eqref{EquationEpsOverThree}}{>} \Big(a-\dfrac{\varepsilon}{3}\Big)^{Nq}$
& $\stackrel{\eqref{EquationEpsOverThree}}{>} \Big(a-\dfrac{2 \varepsilon}{3} \Big)^{(N+l(N))q}$\\[2ex]
& $= \Big(a-\dfrac{2 \varepsilon}{3} \Big)^{n-r}$
& $> \Big(a-\dfrac{2 \varepsilon}{3} \Big)^{n-N-l(N)}$
& $\stackrel{\eqref{EquationNandlN}}{>}(a-\varepsilon)^n$.
\end{tabular}
\end{center}
Thus we see that $a_n > a - \varepsilon$ for almost every $n \in \N$ and the claim follows since $\varepsilon$ was chosen arbitrarily small.
\end{proof}

An application of Proposition~\ref{PropositionCalculusOptimized} in the situation of Lemma~\ref{LemmaSubexpUndistorted} equips us with the following existence criterion for the relative growth rate.

%

\begin{corollary}\label{CorollarySubexpUndistorted}
Let $G$ be a finitely generated group and let $H$ be a finitely generated, subexponentially distorted subgroup of $G$.
If $H$ possesses a polynomially bounded ambiguity function, then the relative growth rate of $H$ exists with respect to every finite generating set of $G$.
\end{corollary}

It would be interesting to know whether there is a finitely generated group $G$ and an undistorted subgroup $H$ of $G$ such that the relative growth rate of $H$ in $G$ does not exist with respect to some finite generating set of $G$.
The following result, which is a direct application of Proposition~\ref{PropositionCalculusOptimized} in the situation of Lemma~\ref{LemmaNewSubexponential}, gives an obstruction to this case.

\begin{corollary}\label{CorollaryNewSubexponential}
Let $G$ be a finitely generated group and let $H$ be a finitely generated undistorted subgroup of $G$.
If $H$ possesses a subexponential ambiguity function, then the relative growth rate of $H$ exists with respect to every finite generating set of $G$.
\end{corollary}

We end this section with the following observation, which was our original motivation for introducing the relative ambiguity function.

\begin{corollary}\label{CorollaryExistSubexpRelAmbiguity}
Let $G$ be a finitely generated group and let $H \leq G$ be an arbitrary subgroup.
If $H$ possesses a subexponential relative ambiguity function in $G$, then the relative growth rate of $H$ exists with respect to every finite generating set of $G$.
\end{corollary}
\begin{proof}
Let $\delta$ be a subexponential relative ambiguity function of $H$ in $G$ with respect to some finite generating set $X$ of $G$.
Recall that this means that there is a constant $\varepsilon \in \N$ and a concatenation map
$\Phi \colon H \times H \rightarrow H,\ (h,k) \mapsto h x_{h,k} k$
with $\abs{x_{h,k}}_X \leq \varepsilon$ for all $h,k \in H$ such that
\[
\abs{B^X_H(m) \times B^X_H(n) \cap \Phi^{-1}(g)} \leq \delta(n)
\]
for all $m,n \in \N$ and every $g \in H$.
Since $\Phi$ restricts to a map
\[
B^X_H(m) \times B^X_H(n) \rightarrow B^X_H(m+n+\varepsilon),
\]
it follows from Dirichlet's box principle that
\[
\beta^{X}_H(m) \beta^{X}_H(n) \leq \delta(n) \cdot \beta^{X}_H(m+n+\varepsilon)
\hspace*{2mm}\text{ for all } m,n \in \N.
\]
Thus the claim follows by applying Proposition~\ref{PropositionCalculusOptimized} on $\beta^{X}_H$.
\end{proof}

\section{Hyperbolic spaces}\label{SectionHyperbolicSpaces}

In this section we recall and introduce some properties of hyperbolic metric spaces.
We start by defining the Gromov product, which can be thought of as the length of a maximal common initial subpath of two geodesics emanating from the same point.

\begin{definition}
Let $(\Omega,d)$ be a metric space and let $o,x,y \in \Omega$ be arbitrary points.
The \emph{Gromov product} of $x$ and $y$ at $o$ is defined by
\[
(x.y)_o = \frac{1}{2}(d(x,o)+d(y,o)-d(x,y)).
\]
\end{definition}

\begin{definition}[Gromov]\label{DefinitionGromovCondition}
A metric space $(\Omega,d)$ is called \emph{$\delta$-hyperbolic} for some constant $\delta \geq 0$ if all points $o,x,y,z \in \Omega$ satisfy
\begin{equation}\label{EquationGromovCondition}
(x.y)_o \geq \min((x.z)_o , (z.y)_o) - \delta.
\end{equation}
In the following we will refer to \eqref{EquationGromovCondition} as the \emph{Gromov condition}.
\end{definition}

Note that the metric space $(\Omega,d)$ in Definition~\ref{DefinitionGromovCondition} is not required to be path-connected.
Nevertheless it will be often useful to require that hyperbolic spaces are geodesic.
Recall that a metric space $(\Omega,d)$ is called \emph{geodesic} if for every two points $x,y \in \Omega$ there is an interval $[a,b] \subset \R$ and an isometric embedding $\gamma \colon [a,b] \rightarrow \Omega$ with $\gamma(a) = x$ and $\gamma(b) = y$.

\begin{definition}\label{DefinitionHypSpGromov}
A geodesic metric space $(\Omega,d)$ is called \emph{hyperbolic} if it is $\delta$-hyperbolic for some $\delta \geq 0$.
\end{definition}

\begin{definition}
Let $(\Omega,d)$ be a metric space and let $I \subset \R$ be an (not necessarily bounded) interval.
Let $\lambda \geq 1$ and $\varepsilon \geq 0$ be some constants.
A map $\gamma \colon I \rightarrow \Omega$ is a \emph{$(\lambda,\varepsilon)$-quasigeodesic} if
\[
\frac{1}{\lambda} \abs{t'-t} - \varepsilon \leq d(\gamma(t'),\gamma(t)) \leq \lambda \abs{t'-t} + \varepsilon
\]
for every $t',t \in I$.
We say that $\gamma$ is a \emph{quasigeodesic} if it is a $(\lambda,\varepsilon)$-quasigeodesic for some $\lambda \geq 1$ and $\varepsilon \geq 0$.
We can replace $I$ by $I \cap \mathbb{Z}$ to get a discrete version of a $(\lambda,\varepsilon)$-quasigeodesic.
\end{definition}

Note that a geodesic is a $(1,0)$-quasigeodesic.
Although geodesics in $\delta$-hyperbolic spaces are generally not determined by their endpoints, we will denote a geodesic connecting two points $x,y$ of a hyperbolic space $(\Omega,d)$ by $[x,y]$.
This should not lead to confusions since, as we will see below, geodesics
sharing the same endpoints stay close together.
Further it will be convenient to identify $[x,y]$ with its image in $(\Omega,d)$, which we will refer to as a \emph{geodesic segment}.

\begin{definition}
Let $(\Omega,d)$ be a metric space.
For any subset $X \subset \Omega$ and any number $\varepsilon > 0$ we denote by $X_{\varepsilon}$ the closed $\varepsilon$-neighbourhood of $X$ in $\Omega$.
Given two subsets $A,B \subset \Omega$, the \emph{Hausdorff distance} between them is defined by
\[
\Hd(A,B) \defeq \inf \{\varepsilon \in \mathbb{R}_{\geq 0} \ \vert \ A \subset B_{\varepsilon},\ B \subset A_{\varepsilon}\}.
\]
\end{definition}

\begin{theorem}[{\cite[Theorem~III.H.1.7]{BridsonHaefliger99}\label{TheoremHausdorffdistance}}]
For all $\delta > 0$, $\lambda \geq 1$, and $\varepsilon \geq 0$ there is a constant $K = K(\delta,\lambda,\varepsilon)$
with the following property:
Let $(\Omega,d)$ be a $\delta$-hyperbolic space, let $\gamma$ be a finite $(\lambda, \varepsilon)$-quasigeodesic in $\Omega$, and let $[x,y]$ be a geodesic segment connecting the endpoints of $\gamma$.
Then the Hausdorff distance $\Hd([x,y],\im(\gamma))$ between every geodesic segment $[x,y]$ and the image $\im(\gamma)$ of $\gamma$ is bounded by $K$.
\end{theorem}

The following lemma tells us that two geodesics emanating from the same point are running parallel to each other if their endpoints are close to each other.

\begin{lemma}[{\cite[Lemma~III.H.1.15]{BridsonHaefliger99}\label{parallelenlemma}}]
Let $(\Omega,d)$ be a $\delta$-hyperbolic space and let $\gamma,\gamma'$ be two geodesics in $\Omega$ emanating from the same point.
If $d(\gamma(t_0),\im(\gamma')) \leq \varepsilon$ for some $\varepsilon \geq 0$ and some $t_0 \geq 0$, then
$d(\gamma(t),\gamma'(t)) \leq 2 \delta \text{ for every } t \leq t_0 - \varepsilon - \delta$.
\end{lemma}

The next definition is attributed to Rips.
It provides us with another important characterization of hyperbolic spaces (see e.g.~\cite[Proposition~III.H.1.22]{BridsonHaefliger99}).

\begin{definition}[Rips]\label{DefinitionRipsCondition}
A geodesic metric space $(\Omega,d)$ is called \emph{hyperbolic} if there is a constant $\delta \geq 0$ such that for every geodesic triangle in $(\Omega,d)$ each of its segments lies in the $\delta$-neighbourhood of the other two segments.
In other words if
\begin{equation}\label{EquationRipsCondition}
[x,y] \subset [x,z]_{\delta} \cup [z,y]_{\delta} \text{ for all } x,y,z \in \Omega
\end{equation}
and every choice of geodesic segments connecting these points.
In the following we will refer to \eqref{EquationRipsCondition} as the \emph{Rips condition}.
\end{definition}

The constant $\delta$ in the Rips condition is not necessarily equal to the constant $\delta$ in the Gromov condition.
In the following we will always assume that $\delta$ is large enough so satisfy both conditions.
To simplify some arguments we will further assume that $\delta$ is a positive integer.

\begin{definition}\label{DefinitionZenter}
Let $(\Omega,d)$ be a $\delta$-hyperbolic space and let $x,y,z \in \Omega$.
For any $\varepsilon \geq 0$ let $C(x,y,z;\varepsilon)$ denote the set of points $p \in \Omega$ such that
\[
p \in [x,y]_{\varepsilon} \cap [y,z]_{\varepsilon} \cap [z,x]_{\varepsilon}
\]
for every choice of geodesic segments $[x,y], [y,z]$ and $[z,x]$.
In the following we will refer to $C(x,y,z;\varepsilon)$ as the \emph{$\varepsilon$-center} of $x,y,z$.
\end{definition}

The following lemma summarizes some properties of $\varepsilon$-centers.

\begin{lemma}\label{LemmaPropertiesOfCenter}
Let $(\Omega,d)$ be a $\delta$-hyperbolic space, let $x,y,z \in \Omega$, and let $\varepsilon \geq 0$.
\begin{enumerate}
\item[\rm{1)}] The set $C(x,y,z;3\delta)$ is non empty.\label{EnumPropertiesOfCenter1}
\item[\rm{2)}] For every $p \in C(x,y,z;\varepsilon)$ we have
\[
(x.y)_z - \varepsilon \leq d(z,p) \leq (x.y)_z + 2 \varepsilon.
\]
\item[\rm{3)}] The diameter of $C(x,y,z;\varepsilon)$ is bounded from above by $7 \varepsilon$.
\item[\rm{4)}] For every choice of geodesic $[z,x]$ there is a point $r \in [z,x]$ with
\[
(x.y)_z - 2\varepsilon \leq d(z,r) \leq (x.y)_z + 3 \varepsilon \text{ for } i = 1,2
\]
that satisfies
\[
d(r,p) \leq 6 \varepsilon \text{ for every } p \in C(x,y,z;\varepsilon).
\]
In particular we have $r \in C(x,y,z;7 \varepsilon)$.
\end{enumerate}
\end{lemma}
\begin{proof}
Let $[x,y], [y,z], [x,z]$ be some choice of geodesic segments and let $t$ be maximal with $[x,y](t) \in [x,z]_{\delta}$.
Then by the Rips condition we have $[x,y](t) \in [y,z]_{\delta}$ and therefore $[x,y](t) \in [x,y]_{\delta} \cap [y,z]_{\delta} \cap [x,z]_{\delta}$.
Now the first claim follows by applying Lemma~\ref{parallelenlemma} to the segments $[x,y], [y,z]$ and $[x,z]$.

To see the second claim we observe that the triangle inequality gives us

%
%

\begin{center}
\begin{tabular}{ll}
$(x.y)_{z}$
& $ = \frac{1}{2}(d(x,z) + d(y,z) - d(x,y))$\\[2ex]
& $\leq \frac{1}{2}(d(x,p) + d(p,z) + d(y,p) + d(p,z) - d(x,p) - d(y,p) + 2 \varepsilon)$\\[2ex]
& $= d(p,z) + \varepsilon.$
\end{tabular}
\end{center}
On the other hand the same argument gives us
\begin{center}
\begin{tabular}{ll}
$d(p,z)$
& $= \frac{1}{2}(d(x,p) + d(p,z) + d(y,p) + d(p,z) - d(x,p) - d(y,p))$\\[2ex]
& $\leq \frac{1}{2}(d(x,z) + 2 \varepsilon + d(y,z) + 2 \varepsilon - d(x,y))$\\[2ex]
& $= (x.y)_{z} + 2 \varepsilon.$
\end{tabular}
\end{center}

For the third claim suppose that $p_1,p_2 \in C(x,y,z;\varepsilon)$.
By the second claim we have
\[
(x.y)_z - \varepsilon \leq d(z,p_i) \leq (x.y)_z + 2 \varepsilon \text{ for } i = 1,2
\]
and therefore $\abs{d(z,p_1) - d(z,p_2)} \leq 3 \varepsilon$.
Let $q_1,q_2 \in [z,x]$ be such that $d(p_1,q_1) \leq \varepsilon$ and $(p_2,q_2) \leq \varepsilon$.
By the triangle inequality we have
\begin{center}
\begin{tabular}{lll}
$(x.y)_z - 2\varepsilon$ 
& $\leq d(z,p_i) - \varepsilon$
& $\leq d(z,p_i) - d(p_i,q_i)$ \\[2ex]
\hspace{-5mm} $\leq d(z,q_i)$
& $\leq d(z,p_i) + d(p_i,q_i)$
& $\leq (x.y)_z + 3 \varepsilon$
\end{tabular}
\end{center}
for $i = 1,2$.
In particular we have
\begin{equation}\label{EquationDistCenterToGrDistPt}
(x.y)_z - 2\varepsilon \leq d(z,q_i) \leq (x.y)_z + 3 \varepsilon \text{ for } i = 1,2
\end{equation}
and therefore $\abs{d(z,q_1) - d(z,q_2)} \leq 5 \varepsilon$.
Since $q_1,q_2$ lie on the same geodesic $[x,z]$, we have $d(q_1,q_2) \leq 5 \varepsilon$.
Now the claim follows from
\[
d(p_1,p_2) \leq d(p_1,q_1) + d(q_1,q_2) + d(q_2,p_2) \leq \varepsilon + 5 \varepsilon + \varepsilon = 7 \varepsilon.
\]
Since $p_1$ and $p_2$ where chosen arbitrarily, the last claim follows from the observation that $q_1$ and $q_2$ satisfy the condition for $r$.
\end{proof}



\begin{lemma}\label{LemmaClosePass}
Let $(\Omega,d)$ be a $\delta$-hyperbolic space and let $w,x,y,z \in \Omega$ be such that
\[
d(x,y) \geq (w.y)_x + (x.z)_y + 41 \delta.
\]
Then every geodesic $[w,z]$ passes through the $((w.y)_x + 31 \delta)$-neighbourhood of $x$.
\end{lemma}
\begin{proof}
By the first claim of Lemma~\ref{LemmaPropertiesOfCenter} we can choose points
\[
p_1 \in C(w,x,y;3\delta) \text{ and } p_2 \in C(x,y,z;3\delta).
\]
Therefore the triangle inequality gives us
\[
d(x,p_i) + d(p_i,y) - 6 \delta \leq d(x,y) \leq d(x,p_i) + d(p_i,y) \text{ for } i = 1,2.
\]
On the other hand the second claim of Lemma~\ref{LemmaPropertiesOfCenter} says that
\[
(w.y)_x - 3\delta \leq d(x,p_1) \leq (w.y)_x + 6\delta
\]
and
\[
(x.z)_y - 3\delta \leq d(y,p_2) \leq (x.z)_y + 6\delta.
\]
Thus from our assumption on $d(x,y)$ we obtain
\begin{center}
\begin{tabular}{lll}
$d(x,p_i) + d(p_i,y)$
& $ \geq d(x,y)$\\[2ex]
\hspace{-5mm} $\geq (w.y)_x + (x.z)_y + 41 \delta$ & $ \geq d(x,p_1) - 6\delta + d(y,p_2) - 6\delta + 41 \delta$\\[2ex]
\hspace{-5mm} $= d(x,p_1) + d(y,p_2) + 29 \delta$\\[2ex]
\end{tabular}
\end{center}
for $i = 1,2$.
For $i = 1$, respectively $i = 2$, we deduce that
\[
d(p_1,y) \geq d(y,p_2) + 29 \delta \text{ and } d(x,p_2) \geq d(x,p_1) + 29 \delta.
\]
By summing up these inequalities, it follows that
\begin{equation}\label{EquationLemmaClosePass1}
d(p_1,y) + d(x,p_2) \geq d(y,p_2) + d(x,p_1) + 58 \delta.
\end{equation}
Recall that by the definition of $C(w,x,y;3\delta)$ we have $p_1 \in [w,x]_{3 \delta} \cap [w,y]_{3 \delta}$.
On the other hand the Rips condition gives us
\[
[w,x] \subset [w,z]_{\delta} \cup [z,x]_{\delta} \text{ and } [w,y] \subset [w,z]_{\delta} \cup [z,y]_{\delta},
\]
so that
\[
p_1 \in ([w,z]_{4 \delta} \cup [z,x]_{4 \delta}) \cap ([w,z]_{4 \delta} \cup [z,y]_{4 \delta}).
\]
Suppose that $p \notin [w,z]_{4 \delta}$.
Then $p_1 \in [z,x]_{4 \delta} \cap [z,y]_{4 \delta}$ and therefore $p_1 \in C(x,y,z;4\delta)$.
The third claim of Lemma~\ref{LemmaPropertiesOfCenter} therefore implies that $d(p_1,p_2) \leq 28 \delta$.
Together with~\eqref{EquationLemmaClosePass1} this shows that the assumption $p \notin [w,z]_{4 \delta}$ leads to a contradiction:
\begin{center}
\begin{tabular}{lll}
$d(p_1,y) + d(x,p_1) + 28 \delta$
& $\geq d(p_1,y) + d(x,p_2)$ \\[2ex]
\hspace{-5mm} $\geq d(y,p_2) + d(x,p_1) + 58 \delta$
& $\geq d(y,p_1) + d(x,p_1) + 58 \delta - 28 \delta$ \\[2ex]
\hspace{-5mm} $= d(y,p_1) + d(x,p_1) + 30 \delta$.
\end{tabular}
\end{center}
Thus we have $p_1 \in [w,z]_{4 \delta}$.
On the other hand the last claim of Lemma~\ref{LemmaPropertiesOfCenter} tells us that there is a point $r \in [x,w]$ with
\[
(w.x)_y - 6 \delta \leq d(x,r) \leq (w.x)_y + 9 \delta \text{ for } i = 1,2
\]
that satisfies $d(r,p_1) \leq 18 \delta$.
Together this implies that $[w,z]$ passes through the $((w.x)_y + (4+18+9)\delta)$-neighbourhood of $x$.
\end{proof}

\section{Acylindrically hyperbolic groups}\label{SectionAcylindricallyHypGrps}


The groups we are mainly interested in are acylindrically hyperbolic groups and products of them.
These groups are known to contain many important and diverse classes of groups such as mapping class groups $\MCG(\Sigma_g)$ of compact orientable surfaces $\Sigma_g$ of genus $g \geq 1$ (see~\cite{Bowditch08} together with~\cite{MasurMinsky99}), outer automorphism groups $\Out(F_r)$ of free groups of rank $r \geq 2$ (see~\cite{BestvinaFeighn10} together with~\cite{Osin16}), and finitely presented groups of deficiency at least $2$ (see~\cite{Osin15}).
Yet acylindrically hyperbolic groups are still accessible for many methods from the theory of hyperbolic groups.
In this section we collect some
results on acylindrically hyperbolic groups.
All the details can be found in~\cite{Osin16}.
We start by defining acylindrical actions on metric spaces.

 
\begin{definition}\label{Def:AcylAction}
Let $G$ be a group that acts
on a metric space $(\Omega,d)$.
The action is called \emph{acylindrical} if for every $\varepsilon > 0$ there are
constants $R,N > 0$ such that for every two points $x,y \in \Omega$ with
$d(x,y) \geq R$, there are at most $N$ elements $g \in G$ such that
$d(x,gx) \leq \varepsilon$ and $d(y,gy) \leq \varepsilon$.
\end{definition}


\begin{definition}
Let $G$ be a group that acts on a hyperbolic space $(\Omega,d)$.
An element $g \in G$ is called \vspace{1mm}
\begin{enumerate}
\item[$\bullet$] \emph{elliptic}, if the set $\{ g^{n}x \ \vert \ n \in \Z\}$
is bounded for some (equivalently, any) $x \in \Omega$.\vspace{2mm}
\item[$\bullet$] \emph{loxodromic}, if the map
$\Z \rightarrow \Omega, \ n \mapsto g^{n}x$ is a quasigeodesic
for some (equivalently, any) $x \in \Omega$.
\end{enumerate}
\end{definition}

Recall that two loxodromic elements are called \emph{independent} if their
fixed point sets in the Gromov boundary $\partial \Omega$ of $\Omega$ are disjoint.
In~\cite[Theorem~1.1]{Osin16} Osin classified acylindrical actions on hyperbolic spaces as follows:

\begin{theorem}\label{TheoremClassificationAcylindric}
Let $G$ be a group that acts acylindrically on a hyperbolic space $(\Omega,d)$.
Then exactly one of the following statements is true:\vspace{1mm}
\begin{enumerate}
\item[\rm{1)}] $G$ has bounded orbits.
\item[\rm{2)}] $G$ is virtually cyclic and contains a loxodromic element.
\item[\rm{3)}] $G$ contains infinitely many independent loxodromic elements.
\end{enumerate}
\end{theorem}

The first two cases in Theorem~\ref{TheoremClassificationAcylindric} are referred to as \emph{elementary actions} and
are excluded from providing examples of acylindrically hyperbolic groups.

\begin{definition}\label{Def:AcylHypGrp}
A group $G$ is called \emph{acylindrically hyperbolic} if $G$ admits a non-elementary acylindrical action on a hyperbolic space.
\end{definition}

The following theorem allows us to switch from acylindrical actions on general hyperbolic spaces to acylindrical actions on hyperbolic Cayley graphs.
This will be crucial for us since general hyperbolic spaces are not symmetric enough for our purposes.

\begin{theorem}[{\cite[Theorem~1.4]{Osin16}\label{TheoremAcylActionOnCayleyGraph}}]
Let $G$ be a group and let $g \in G$ be an arbitrary element.
The following conditions are equivalent.
\begin{enumerate}
\item[\rm{1)}] There exists a (possibly infinite) generating set $X$ of $G$ such that the corresponding Cayley
graph $\Gamma(G, X)$ is hyperbolic, the natural action of $G$ on $\Gamma(G, X)$ is
acylindrical, and $g$ is loxodromic.
\item[\rm{2)}] There exists an acylindrical action of $G$
on a hyperbolic space $(\Omega,d)$ such that $g$ is loxodromic.
\end{enumerate}
\end{theorem}

Note that Theorem~\ref{TheoremAcylActionOnCayleyGraph} has the remarkable consequence that an acylindrical action on a hyperbolic space can be replaced by a cocompact one.
If the group $G$ in Theorem~\ref{TheoremAcylActionOnCayleyGraph} is finitely generated it is sometimes useful to assume that $X$ contains a finite generating set of $G$.
%
The following lemma implies in particular that we can add a finite generating set to $X$ without violating
the acylindricity of the action of $G$ on $\Gamma(G, X)$.

\begin{lemma}[{\cite[Lemma~5.1]{Osin16}\label{LemmaChangeOfGenratingSet}}]
For any group $G$ and any generating sets $X$ and $Y$ of $G$ such that
\[
\sup \limits_{x \in X} \abs{x}_Y < \infty \text{ and } \sup \limits_{y \in Y} \abs{y}_X < \infty,
\]
the following hold.
\begin{enumerate}
\item[\rm{1)}] $\Gamma(G, X)$ is hyperbolic if and only if $\Gamma(G, Y)$ is hyperbolic.
\item[\rm{2)}] $G$ acts acylindrically on $\Gamma(G, X)$ if and only if $G$ acts acylindrically on $\Gamma(G, Y)$.
\end{enumerate}
\end{lemma}

In general it may happen
that the element $g$ in Theorem~\ref{TheoremAcylActionOnCayleyGraph} is loxodromic for some but not for every acylindrical action on a hyperbolic space.
For our purposes it will be sufficient to know that an element is loxodromic for some action on a hyperbolic space.

\begin{definition}[{\cite[Definition~6.4]{Osin16}\label{DefinitionGeneralizedLoxo}}]
Let $G$ be an arbitrary group.
An element $g \in G$ is called \emph{generalized loxodromic} if there is an acylindrical action of $G$ on a hyperbolic space $(\Omega,d)$ such that $g \in G$ is a loxodromic element for this action.
\end{definition}

\section{Existence of the growth rate}\label{SectionExistence}

The aim of this section is to prove Theorem~\ref{IntroMainTheorem} from the introduction.
As a consequence we will answer Question~\ref{QuestionExistenceForHyperbolic} affirmatively.
We start with the following observation.

\begin{lemma}\label{LemmaAndDefConcatenation}
Let $G$ be a group with a generating set $Y$ such that $\Gamma(G,Y)$ is $\delta$-hyperbolic and the canonical action of $G$ on $\Gamma(G,Y)$ is acylindrical and non-elementary.
Then for every two independent loxodromic elements $g$ and $h$ there is a constant $A \geq 0$ such that for every $n \in \N$ and every pair of elements $(u,v) \in G \times G$ there is at least one element $x_{u,v} \in \{ g^{n},g^{-n},h^{n},h^{-n} \}$ that satisfies
\[
(u^{-1}.x_{u,v})_1 \leq A \text{ and } (v.x_{u,v}^{-1})_1 \leq A.
\]
\end{lemma}
\begin{proof}
Since $g$ and $h$ are independent, there is a constant $C \geq 0$ such that the Gromov product $(x.y)_1$ of two distinct elements $x,y \in \{ g^{n},g^{-n},h^{n},h^{-n} \}$ with respect to the identity is bounded above by $C$.
From the Gromov condition we therefore obtain
\begin{equation}\label{Eq:GromovCross}
C \geq (x.y)_{1} \geq \min \{ (x.w)_1 , (w.y)_1 \} - \delta
\end{equation}
for every $w \in G$ and every $n \in \N$.
Note that if $C < (g^{n}.w)_1 - \delta$ or $C < (g^{-n}.w)_1 - \delta$, then~\eqref{Eq:GromovCross} implies that $C \geq (h^n.w)_1 - \delta$ and $C \geq (h^{-n}.w)_1 - \delta$.
Hence we have
\begin{equation}\label{Eq:GromovCross2}
C \geq \min \{ \max \{ (g^n.w)_1 , (w.g^{-n})_1 \}, \max \{ (h^n.w)_1 , (w.h^{-n})_1 \} \} - \delta
\end{equation}
for every $w \in G$ and every $n \in \N$.
Let $u,v \in G$ be as in the lemma.
Then, on the one hand,~\eqref{Eq:GromovCross2} provides us with an element
$x \in \{ g^{n},g^{-n},h^{n},h^{-n} \}$ such that
\[
C \geq (x.u^{-1})_1 - \delta \text{ and } C \geq (x^{-1}.u^{-1})_1 - \delta
\]
and on the other hand~\eqref{Eq:GromovCross} tells us that $C \geq \min \{(x.v)_1 ,(x^{-1}.v)_1\} - \delta$.
If $C \geq (x^{-1}.v)_1 - \delta$, then the claim follows by setting $A \defeq C + \delta$ and $x_{u,v} \defeq x$.
If $C \geq (x.v)_1 - \delta$ we can replace $x$ by $x^{-1}$ and proceed as before.
\end{proof}

Since the elements $g,h \in G$ in Lemma~\ref{LemmaAndDefConcatenation} are loxodromic, every $R \geq 0$ gets exceeded by $\abs{g^n}_Y$ and $\abs{h^n}_Y$ for some appropriate $n \in \N$.
Suppose that $G$ is also finitely generated with a finite generating set $X$.
Then we obtain the following corollary by setting $\varepsilon \defeq \max \{ \abs{g^n}_X,\abs{h^n}_X \}$ and by choosing $B^{X}_{H}(\varepsilon)$ to be a set of connecting pieces.

\begin{corollary}\label{CorollaryExistNiceAmbiguityFct}
Let $G$ be a finitely generated group with a finite generating set $X$ and let $H$ be a subgroup of $G$.
Suppose that $G$ has a generating set $Y$ such that $\Gamma(G,Y)$ is hyperbolic, the canonical action of $G$ on $\Gamma(G,Y)$ is acylindrical, and that $H$ contains two independent loxodromic elements.
Then there is a constant $A \geq 0$ such that for every $R \geq 0$ there is
a concatenation map
\[
\Phi \colon H \times H \rightarrow H,\ (h,k) \mapsto h x_{h,k} k,
\]
with the following properties.
For all $h,k \in H$ the connecting piece $x_{h,k}$ has length $\abs{x_{h,k}}_{Y} \geq R$
and satisfies
$(h^{-1}.x_{h,k})_1, (k.x_{h,k}^{-1})_1 \leq A$.
\end{corollary}


It will be useful to observe that the Gromov product is invariant under isometry.

\begin{remark}\label{RemarkLeftInvGromovProd}
Let $G$ be a group that acts by isometries on a metric space $(\Omega,d)$.
Then the Gromov product satisfies $(x.y)_z = (gx.gy)_{gz}$ for every $g \in G$ and all points $x,y,z \in \Omega$.
\end{remark}

\begin{proposition}\label{PropositionCloseMiddleCross}
Let $G$ be a finitely generated group with a finite generating set $X$ and let $H$ be a subgroup of $G$.
Suppose that $G$ has a generating set $Y$ such that $\Gamma(G,Y)$ is hyperbolic, the canonical action of $G$ on $\Gamma(G,Y)$ is acylindrical, and that $H$ contains two independent loxodromic elements.
Then there is a constant $C \geq 0$ such that for every $R \geq 0$ there is a concatenation map
\[
\Phi \colon H \times H \rightarrow H,\ (h,k) \mapsto h x_{h,k} k
\]
with the following properties.
For all $h,k \in H$ the connecting piece $x_{h,k}$ has length $\abs{x_{h,k}}_{Y} \geq R$ and every geodesic segment $[1,\Phi(h,k)]$ in $\Gamma(G,Y)$
passes through the $C$-neighbourhoods of $h$ and $h x_{h,k}$.
\end{proposition}
\begin{proof}
From Corollary~\ref{CorollaryExistNiceAmbiguityFct} we know that there is a constant $A \geq 0$ and a concatenation map
\[
\Phi \colon H \times H \rightarrow H,\ (h,k) \mapsto h x_{h,k} k,
\]
that satisfies
\[
\abs{x_{h,k}}_Y \geq \max \{2A + 41 \delta, R \} \text{ and } (h^{-1}.x_{h,k})_1, (k.x_{h,k}^{-1})_1 \leq A \text{ for all } h,k \in H.
\]
We want to apply Lemma~\ref{LemmaClosePass} on the points $1,h,hx_{h,k}$ and $h x_{h,k} k$.
Thus we have to check that
\[
d_Y(h,hx_{h,k}) \geq (1.hx_{h,k})_h + (h.h x_{h,k} k)_{hx_{h,k}} + 41 \delta
\]
in order to prove that every geodesic $[1,h x_{h,k} k]$ passes through the $((1.ux_{h,k})_h + 31 \delta)$-neighbourhood of $h$.
By Remark~\ref{RemarkLeftInvGromovProd} we have
\[
(1.h x_{h,k})_h = (h^{-1}.x_{h,k})_1 \leq A
\]
and similarly
\begin{center}
\begin{tabular}{ll}
$(h.hx_{h,k}k)_{hx_{h,k}}$
& $= ((hx_{h,k})^{-1} h.(hx_{h,k})^{-1} hx_{h,k}k)_{(hx_{h,k})^{-1} hx_{h,k}}$\\[2ex]
& $= (x_{h,k}^{-1}.k)_{1} \leq A$.
\end{tabular}
\end{center}
Thus we obtain
\[
d_Y(u,u x_{u,v}) = \abs{x_{u,v}}_Y \geq 2A + 41 \delta
\]
and by setting $C \defeq A + 31 \delta$ it follows that every geodesic $[1,h x_{h,k} k]$ passes through the $C$-neighbourhood of $h$.
Symmetric argumentation shows that every geodesic $[h x_{h,k} k,1]$ passes through the $C$-neighbourhood of $h x_{h,k}$.
\end{proof}

\begin{proposition}\label{PropositionBoundedPreimage}
Let $G$ be a finitely generated group with a finite generating set $X$ and a subgroup $H$.
Suppose that $G$ has a generating set $Y$ such that $\Gamma(G,Y)$ is hyperbolic, the canonical action of $G$ on $\Gamma(G,Y)$ is acylindrical, and that $H$ contains two independent loxodromic elements.
Then $H$ possesses a linearly bounded relative ambiguity function in $G$.
\end{proposition}
\begin{proof}
Due to Lemma~\ref{LemmaChangeOfGenratingSet} we can assume that $X$ is a subset of $Y$.
Further Proposition~\ref{PropositionCloseMiddleCross} provides us with a constant $C \geq 0$ such that for every $R \geq 0$ there is a concatenation map
\[
\Phi^{(R)} \colon H \times H \rightarrow H,\ (h,k) \mapsto h x^{(R)}_{h,k} k
\]
with respect to a set $B^X_H(\varepsilon^{(R)})$ of connecting pieces with the following properties.
For all $h,k \in H$ the connecting piece $x^{(R)}_{h,k}$ has length $\abs{x^{(R)}_{h,k}}_{Y} \geq R$ and every geodesic segment $[1,\Phi^{(R)}(h,k)]$
passes through the $C$-neighbourhoods of $h$ and $h x^{(R)}_{h,k}$ in $\Gamma(G,Y)$.
Since the action of $G$ on $\Gamma(G,Y)$ is acylindrical, we can choose positive constants $R_0$ and $N_0$ such that for every two points $x,y \in \Gamma(G,Y)$ with $d_Y(x,y) \geq R_0$ there are less than $N_0$ elements $g \in G$ that satisfy
\[
d_Y(x,g x) \leq 6 C \text{ and } d_Y(y,g y) \leq 6 C.
\]
We fix the concatenation map $\Phi \defeq \Phi^{(R_0)}$ and we set $x_{h,k} \defeq x^{(R_0)}_{h,k}$ and $\varepsilon \defeq \varepsilon^{(R_0)}$.
Let $w \in H$ and let $m,n \in \N$.
We have to construct an upper bound for the number of pairs $(h,k) \in B^X_H(m) \times B^X_H(n)$ that satisfy $h x_{h,k} k = w$.
Note that for a fixed $h \in H$ there are at most $D \defeq \abs{B^X_H(\varepsilon)}$ elements $k \in H$ with $h x_{h,k} k = w$.
Thus it suffices to estimate the number of elements $h \in B^X_H(m)$ with $h x_{h,k} k = w$ for some $k \in B^X_H(n)$.
For this we fix a geodesic $[1,w]$ in $\Gamma(G,Y)$.
Then for every pair $(h,k) \in B^X_H(m) \times B^X_H(n)$ with $h x_{h,k} k = w$ there is an element $g \in G$ that that lies on $[1,w] = [1,h x_{h,k} k]$ and satisfies $d_Y(h,g) \leq C$.
We therefore obtain
\begin{align*}
d_Y(g,w) & \leq d_Y(g , h) + d_Y(h , h x_{h,k} k)\\
& \leq C + d_Y(1 , x_{h,k} k)\\
& \leq C + \abs{x_{h,k}}_Y + \abs{k}_Y\\
& \leq C + \abs{x_{h,k}}_X + \abs{k}_X\\
& \leq C + \varepsilon + n.
\end{align*}
In particular we see that there are at most $C + \varepsilon + n$ elements $g \in G \cap [1,w]$ for which there is a pair $(h,k) \in B^X_H(m) \times B^X_H(n)$ with $d_Y(h,g) \leq C$ and $h x_{h,k} k = w$.
We fix such an element $g_0 \in G$.
Suppose that there are $D \cdot N_0$ elements $h \in B^X_H(m) \cap B^Y_H(g_0,C)$ with $\Phi(h,k) = w$ for some $k \in B^X_H(n)$.
Then by Dirichlet's box principle there exists an element $x \in B^X_H(\varepsilon)$ such that there are $N_0$ pairs
\[
(h_1,k_1), \ldots, (h_{N_0},k_{N_0}) \in (B^X_H(m) \cap B^Y_H(g_0,C)) \times B^X_H(n)
\]
that satisfy
\[
h_i x_{h_i,k_i} k_i = h_i x k_i = w \text{ for all } 1 \leq i \leq N_0.
\]
By the choice of $\Phi$ we can fix elements $v_i \in [1,w] = [1,h_i x k_i]$ with $d_Y(v_i,h_i x) \leq C$ for every $1 \leq i \leq N_0$.
We therefore obtain
\[
d_Y(g_0,v_i) \leq d_Y(g_0,h_i) + d_Y(h_i,h_i x) + d_Y(h_i x,v_i) \leq 2C + \abs{x}_Y
\]
and
\[
d_Y(g_0,v_i) \geq d_Y(h_i,h_i x) - d_Y(g_0,h_i) - d_Y(h_i x,v_i) \geq \abs{x}_Y - 2 C.
\]
Thus we have
\[
\abs{x}_Y - 2 C \leq d_Y(g_0,v_i) \leq \abs{x}_Y + 2 C \text{ for all } 1 \leq i \leq N_0.
\]
Since every $v_i$ lies on the same geodesic $[1,w]$, it follows that $d_Y(v_i,v_j) \leq 4 C$ for all $1 \leq i,j \leq N_0$ and therefore
\[
d_Y(h_1 x,h_j x) \leq d_Y(h_1 x,v_1) + d_Y(v_1,v_j) + d_Y(v_j,h_j x)
\leq 6C.
\]
Note that the triangle inequality further gives us
\[
d_Y(h_1,h_j) \leq d_Y(h_1,g_0)+d_Y(g_0,h_j) \leq 2 C.
\]
By setting $\hat{h}_j \defeq h_j^{-1} h_1$, the above inequalities can be rewritten as
\[
d_Y(\hat{h}_j,1) \leq 2 C \text{ and } d_Y(\hat{h}_j x, x) \leq 6C
\]
for every $j \in \{1,\ldots,N_0\}$.
Since $d_Y(1,x) \geq R_0$, this contradicts the acylindricity of the action.
Hence there are less than $D \cdot N_0$ elements $h \in B^X_H(m) \cap B^Y_H(g_0,C)$ with $\Phi(h,k) = w$ for some $k \in B^X_H(n)$.
Since we had at most $C + \varepsilon + n$ choices for $g_0$, we see that there are less than $D \cdot N_0 \cdot (C + \varepsilon + n)$ elements $h \in B^X_H(m)$ that satisfy the equation $\Phi(u,v) = w$ for some appropriate elements $k \in B^X_H(n)$.
As there are at most $D$ possible choices $k \in H$ with $h x_{h,k} k = w$ for some fixed $h \in B^X_H(m)$, we see that the number of all pairs $(h,k) \in B^X_H(m) \times B^X_H(n)$ satisfying $h x_{h,k} k = w$ is bounded above by $D^2 \cdot N_0 \cdot (C + \varepsilon + n)$.
Now the claim follows since
\[
f \colon \N \rightarrow \N,\ n \mapsto D^2 \cdot N_0 \cdot (C + \varepsilon + n)
\]
defines a linearly bounded relative ambiguity function for $H$ in $G$.
\end{proof}

A careful analysis of the proof of Proposition~\ref{PropositionBoundedPreimage} reveals that the concatenation map $\Phi \colon H \times H \rightarrow H,\ (h,k) \mapsto h x_{h,k} k$ that corresponds to the linearly bounded relative ambiguity function of $H$ extends to a map $\widetilde{\Phi} \colon G \times G \rightarrow G,\ (h,k) \mapsto h x_{h,k} k$ with $x_{h,k} \in H$ for all $h,k \in G$.
Indeed, tracing back the argument to Lemma~\ref{LemmaAndDefConcatenation}, we see that $x_{h,k}$ was constructed for any two elements $h,k \in G$ as a power of one of two loxodromic elements that were subsequently chosen from $H$.
It then follows that the ambiguity function of $G$ that corresponds to $\widetilde{\Phi}$ is linearly bounded.
We summarize this as follows.

\begin{remark}\label{RemarkExtensionOfTheRelativeAmbFct}
In the situation of Proposition~\ref{PropositionBoundedPreimage}, the concatenation map $\Phi$ of $H$ that corresponds to the linearly bounded relative ambiguity function of $H$, extends to a concatenation map $\widetilde{\Phi} \colon G \times G \rightarrow G,\ (h,k) \mapsto h x_{h,k} k$ with $x_{h,k} \in H$ for all $h,k \in G$, such that $\widetilde{\Phi}$ induces a linearly bounded ambiguity function of $G$.
\end{remark}

\begin{theorem}\label{TheoremEmbeddingInvariant}
Let $G$ be a finitely generated acylindrically hyperbolic group and let $H$ be a subgroup of $G$.
If $H$ contains a generalized loxodromic element of $G$, then $H$ has a linearly bounded relative ambiguity function in $G$.
\end{theorem}
\begin{proof}
From Theorem~\ref{TheoremAcylActionOnCayleyGraph} we know that $G$ has a generating set $Y$ such that $\Gamma(G,Y)$ is hyperbolic, the canonical action of $G$ on $\Gamma(G,Y)$ is acylindrical, and $h$ is a loxodromic element for this action.
Then Theorem~\ref{TheoremClassificationAcylindric} implies that there are two cases to consider for $H$.
Either $H$ is virtually cyclic or $H$ contains infinitely many independent loxodromic elements.
Recall that in both cases Lemma~\ref{LemmaChangeOfGenratingSet} allows us to assume that $Y$ contains a finite generating set $X$ of $G$.
In the first case $H$ is quasiisometrically embedded in $\Gamma(G,Y)$.
In particular $H$ is quasiisometrically embedded in $\Gamma(G,X)$ and there is a constant $C > 0$ such that $\beta^X_H(n) \leq C \cdot n$ for every $n \in \N$.
Thus we see that the limit $\lim \limits_{n \rightarrow \infty} \sqrt[n]{\beta^X_H(n)}$ exists and is equal to $1$.
The second case is a direct consequence from Proposition~\ref{PropositionBoundedPreimage}.
\end{proof}

Recall from Corollary~\ref{CorollaryExistSubexpRelAmbiguity} that the relative growth rate of $H$ exists with respect to every finite generating set of $G$,
if $H$ has a subexponential relative ambiguity function in $G$.
As a consequence of Theorem~\ref{TheoremEmbeddingInvariant}, we therefore obtain Theorem~\ref{IntroMainTheorem} from the introduction.

\begin{theorem}\label{TheoremMain}
Let $G$ be a finitely generated acylindrically hyperbolic group and let $H$ be a subgroup of $G$.
If $H$ contains a generalized loxodromic element of $G$, then the relative growth rate of $H$ exists with respect to every finite generating set of $G$.
\end{theorem}

\begin{remark}\label{RemarkLoxoIsNess}
The condition in Theorem~\ref{TheoremMain} that $H$ contains a generalized loxodromic element cannot be removed.
To see this let $K$ be a finitely generated group that contains a subgroup $L$ such that $\liminf \limits_{n \rightarrow \infty} \sqrt[n]{\beta^Y_L(n)} \neq \limsup \limits_{n \rightarrow \infty} \sqrt[n]{\beta^Y_L(n)}$ for some finite generating set $Y$ of $K$.
The existence of such groups was proven by Olshanskii in~\cite[Remark~3.1]{Olshanskii17}.
We consider the free product $\widetilde{K} \defeq K \ast \Z$ and the canonical embedding $\widetilde{L}$ of $L$ in $\widetilde{K}$.
The group $\widetilde{K}$ is acylindrically hyperbolic by~\cite[Example 2.10]{Osin16}.
However, if $\Z$ is generated by some element $t$, we see that $\widetilde{Y} \defeq Y \cup \{t\}$ generates $\widetilde{K}$ and that $\beta^{\widetilde{Y}}_{\widetilde{L}}(n) = \beta^{Y}_{L}(n)$ for every $n \in \N$.
Thus we see that the relative growth rate of $\widetilde{L}$ does not exist with respect to the generating set $\widetilde{Y}$ of $\widetilde{K}$.
\end{remark}


We are now ready to answer Question~\ref{QuestionExistenceForHyperbolic} from the introduction.

\begin{theorem}\label{TheoremExistenceForHyperbolicGroups}
Let $G$ be a hyperbolic group and let $H \leq G$ be an arbitrary subgroup.
Then the relative growth rate of $H$ exists with respect to every finite generating set $X$ of $G$.
\end{theorem}

In view of this theorem, it follows that relative growth functions of subgroups of hyperbolic groups do not behave as wild as one might expect from the Rips construction~\cite{Rips82}, which tells us that hyperbolic groups may contain extremely distorted subgroups.

\begin{proof}[Proof of Theorem~\ref{TheoremExistenceForHyperbolicGroups}]
The canonical action of $G$ on its Cayley graph $\Gamma(G, X)$ is proper and hence acylindrical.
Further every element of $G$ of infinite order is acting loxodromically on $\Gamma(G, X)$ (see~\cite[Corollary~III.$\Gamma$.3.10]{BridsonHaefliger99}).
Thus by Theorem~\ref{TheoremMain} the relative growth rate of $H$ exists with respect to $X$ if $H$ contains an element of infinite order.
Suppose that $H$ does not contain any element of infinite order.
In this case Theorem~\ref{TheoremClassificationAcylindric} tells us that the action of $H$ on $\Gamma(G, X)$ has bounded orbits.
Since $X$ is finite, it follows that $H$ is finite so that the relative growth rate of $H$ obviously exists.
\end{proof}

Note that the case $H = G$ in Theorem~\ref{TheoremEmbeddingInvariant} provides us with the following invariant of acylindrically hyperbolic groups.

\begin{corollary}\label{CorollaryAcylHypPossNonRelLinAmb}
Every finitely generated acylindrically hyperbolic group $G$ possesses a linearly bounded relative ambiguity function.
\end{corollary}

\section{The case of products}\label{SectionProducts}

Our next goal is to generalize the results from the previous section to products of acylindrically hyperbolic groups.
This will allow us to study relative growth rates in arbitrary right-angled Artin groups.
We start by introducing some notation that we will use throughout the remaining paper.

\begin{notation}\label{NotationEmbeddingProjection}
Let $I = \{ 1,\ldots,m \}$ be a set of indices.
For each $i \in I$ let $G_i$ be a finitely generated group with a finite generating set $X_i$.
For every subset $J \subset I$ we denote by
\begin{center}
$\iota_J \colon \prod \limits_{j \in J} G_j \rightarrow \prod \limits_{i \in I} G_i
\text{ and }
\psi_J \colon \prod \limits_{i \in I} G_i \rightarrow \prod \limits_{j \in J} G_j$
\end{center}
the canonical embedding resp. the canonical projection.
If $J$ is a singleton we will write $\iota_{i} \defeq \iota_{ \{ i \} }$ and $\psi_{i} \defeq \psi_{ \{ i \} }$.
\end{notation}

We will make frequent use of the following easy observation.

\begin{remark}\label{RemarkUnionOfGenSets}
Let $X = \bigcup \limits_{i \in I} \iota_i(X_i)$ be the union of the embeddings of the generating sets in Notation~\ref{NotationEmbeddingProjection}.
Then the length of an element $(g_i)_{i \in I} \in \prod \limits_{i \in I} G_i$ with respect to $d_X$ is given by $\abs{(g_i)_{i \in I}}_X = \sum \limits_{i \in I}^{m} \abs{g_i}_{X_i}$.
\end{remark}

\begin{lemma}\label{LemmaSmallKernel}
Let $G$ be a finitely generated group with a finite generating set $X$.
Let $H,K \leq G$ be subgroups such that there is an epimorphism $\psi \colon H \rightarrow K$ that satisfies $\abs{\psi(h)}_X \leq \abs{h}_X$ for every $h \in H$.
Suppose that $f_K \colon \N \rightarrow \N$ is a relative ambiguity function for $K$ with respect to $X$.
Then, up to linear equivalence, $H$ has a relative ambiguity function $f_H$ with respect to $X$ that is bounded above by $f_K \cdot \beta_{\ker(\psi)}^{X}$.
\end{lemma}
\begin{proof}
By the definition of $f_K$ there is a constant $\varepsilon_K \in \N$ and a concatenation map
$\Phi_K \colon K \times K \rightarrow K,\ (h,k) \mapsto h x_{h,k} k$
with $\abs{x_{h,k}}_X \leq \varepsilon_K$ for all $h,k \in K$ such that
\[
\abs{B^X_K(m) \times B^X_K(n) \cap \Phi_K^{-1}(g)} \leq f_K(n)
\]
for all $m,n \in \N$ and every $g \in K$.
For each $x \in B^X_K(\varepsilon_K)$ we pick a preimage $x'$ of $x$ with respect to $\psi$ and let $\varepsilon_H \defeq \max \limits_{x \in B^X_K(\varepsilon_K)} \abs{x'}_{X}$ denote the maximal length of these preimages.
We claim that the concatenation map
\[
\Phi_H \colon H \times H \rightarrow H,\ (h,k) \mapsto h x'_{\psi(h),\psi(k)} k
\]
gives rise to a relative ambiguity function for $H$ with respect to $X$ that satisfies the requirements.
Since by construction $\abs{x'_{\psi(h),\psi(k)}}_X \leq \varepsilon_H$, it remains to construct an upper bound for a function $f_H \colon \N \rightarrow \N$ that satisfies
\[
\abs{B^X_H(m) \times B^X_H(n) \cap \Phi_H^{-1}(g)} \leq f_H(n)
\]
for all $m,n \in \N$ and every $g \in H$.
To see this we fix two numbers $m,n \in \N$ and an element $g \in B^{X}_{H}(m+n+\varepsilon_H)$.
Let $N \defeq \ker(\psi)$.
We claim that the number of cosets $kN \in H/N$ with $h x'_{\psi(h),\psi(k)} k = g$ for some $h \in B^{X}_{H}(m)$ and a representative $k \in B^{X}_{H}(n)$ is bounded by $\varepsilon_K(n)$.
Indeed, observe that distinct such cosets $k N$ induce distinct ways of writing $\psi(g)$ as
\[
\psi(h) \cdot x_{\psi(h),\psi(k)} \cdot \psi(k) = \psi(h x'_{\psi(h),\psi(k)} k) = \psi(g).
\]
Thus the claim follows from $\psi(B^X_H(n)) \subseteq B^X_K(n)$, which is a direct consequence of our assumption that $\abs{\psi(h)}_X \leq \abs{h}_X$ for every $h \in H$.
Suppose now that $k \in B^{X}_{H}(n)$ is a fixed element with
$h x_{h,k}' k = g$
for some $h \in B^{X}_{H}(m)$.
We want to estimate the cardinality $\abs{kN \cap B^{X}_{H}(n)}$ of the set of elements in $B^{X}_{H}(n)$ that represent $k$ mod $N$.
In order to do so, we observe that
\begin{equation}\label{EquationDifferentRepresentatives}
\abs{kN \cap B^{X}_{H}(n)} = \abs{N \cap k^{-1}B^{X}_{H}(n)} \leq \abs{B^{X}_{N}(2n)} = \beta^{X}_{N}(2n),
\end{equation}
where the last inequality follows from $\abs{k}_{X} \leq n$.
By combining~\eqref{EquationDifferentRepresentatives} with the above observation on the number of cosets $kN \in H/N$ with $h x'_{\psi(h),\psi(k)} k = g$ for some $h \in B^{X}_{H}(m)$, we see that the number of elements $k \in B^{X}_{H}(n)$ satisfying $h x_{h,k}' k$ for some $h \in B^{X}_{H}(m)$ is bounded above by $\varepsilon_K(n) \cdot \beta^{X}_{N}(2n)$.
Since there are at most $\beta^{X}_{H}(\varepsilon_H)$ elements $h \in B^{X}_{H}(m)$ with $h x'_{\psi(h),\psi(k)} k = g$ for some fixed $k$, Dirichlet's box principle implies that the relative ambiguity function $f_H$ for $H$ with respect to $X$ can be defined as
\[
f_H \colon \N \rightarrow \N,\ n \mapsto \beta^{X}_{H}(\varepsilon_H) \cdot \varepsilon_K(n) \cdot \beta^{X}_{N}(2n).
\]
Now the claim follows since, up to linear equivalence, $f_H$ is bounded above by $f_K \cdot \beta_{\ker(\psi)}^{X}$.
\end{proof}


\begin{lemma}\label{LemmaDecomposition}
Let $I = \{ 1,\ldots,m \}$ be a set of indices.
For each $i \in I$ let $G_i$ be a group with a generating set $Y_i$ such that $\Gamma(G_i,Y_i)$ is hyperbolic and $G_i$ acts acylindrically on $\Gamma(G_i,Y_i)$.
Let $H \leq \prod \limits_{i \in I} G_i$ be a subgroup and suppose that $\psi_i(H) \leq G_i$ contains a loxodromic element of $G_i$ for every $i \in I$.
Then there is a decomposition of $I$ into disjoint subsets $I_1$ and $I_2$ such that the following holds.
\begin{enumerate}
\item[\rm{1)}] The kernel of the restricted projection $(\psi_{I_1})_{\vert H} \colon H \rightarrow \prod \limits_{i \in I_1} G_i$ is finite.
\item[\rm{2)}] For every index $i \in I_1$ there is a loxodromic element $g_i \in G_i$ such that $H$ contains an element $(h_j)_{j \in I}$ with $h_i = g_i$ and $h_j = 1$ for $j \in I_1 \setminus \{ i \}$.
\end{enumerate}
\end{lemma}
\begin{proof}
Since $(\psi_{I})_{\vert H}$ is the identity on $H$, we can choose a minimal subset $I_1 \subset I$ such that $(\psi_{I_1})_{\vert H}$ has a finite kernel.
Without loss of generality we may
assume that $I_1$ is a proper subset of $I$.
Note that this implies that the kernel $N_i \defeq \ker((\psi_{I_1 \setminus \{i\}})_{\vert H})$ is infinite for every $i \in I_1$.
Note further that each such kernel $N_i$ embeds as a subgroup of $G_i$ via $\psi_{i}$ and as such acts acylindrically on $\Gamma(G_i,Y_i)$.
We can therefore apply Theorem~\ref{TheoremClassificationAcylindric} which tells us that $\psi_{i}(N_i)$ contains a loxodromic element if and only if it is unbounded in $\Gamma(G_i,Y_i)$.
Suppose the contrary, i.e.\ that $\psi_{i}(N_i) \subset \Gamma(G_i,Y_i)$ is contained in $B_{G_i}^{Y_i}(\varepsilon)$ for some $\varepsilon \geq 0$.
Due to our assumptions there is an element $(h_j)_{j \in I} \in H$ such that $h_i$ is loxodromic.
By conjugating an arbitrary element $(1, \ldots,1,g_i,1, \ldots,1) \in N_i$ with $(h_j^n)_{j \in I}$, we obtain $(1, \ldots,1,h_i^n g_i h_i^{-n},1, \ldots,1) \in N_i$ and we claim that we can choose $g_i$ and $n$ such that $\abs{h_i^{n} h h_i^{-n}}_Y$ is not contained in $B_{G_i}^{Y_i}(\varepsilon)$.
Since the action of $G_i$ on $\Gamma(G_i,Y_i)$ is acylindrical, there are constants $N,R > 0$ such that for every pair of elements $x,y \in G_i$ with $d_Y(x,y) \geq R$ we have
\[
\abs{\{ g \in G_i \ \vert \ d_Y(x,gx) \leq \varepsilon \text{ and } d_Y(y,gy) \leq \varepsilon \}} \leq N.
\]
As $h_i$ is loxodromic, we can choose $n \in \N$ big enough to satisfy $\abs{h_i^{n}}_Y \geq R$.
Then, by setting $x = 1$ and $y = h_i^{n}$, we see that every element $h \in N_i$ satisfies
\[
d_Y(1,h 1) \leq \varepsilon \text{ and } d_Y(h_i^{n},h h_i^{n}) \leq \varepsilon.
\]
But this is a contradiction since $N_i$, and therefore $\psi_i(N_i)$, is infinite.
\end{proof}

\begin{proposition}\label{PropositionTheSimpleCase}
Let $I$ be a finite set of indices.
For each $i \in I$ let $G_i$ be a finitely generated group with a generating set $Y_i$ such that $\Gamma(G_i,Y_i)$ is hyperbolic and $G_i$ acts acylindrically on $\Gamma(G_i,Y_i)$.
Suppose that $G \defeq \prod \limits_{i \in I} G_i$ contains a subgroup $H$ with the following property.
For each $i \in I$ there are elements $g^{(i)} \defeq (g^{(i)}_j)_{j \in I}$ and $h^{(i)} \defeq (h^{(i)}_j)_{j \in I} \in H$ such that $g^{(i)}_{i}$ and $h^{(i)}_{i}$ are independent loxodromic elements in $G_i$ with $g^{(i)}_{j},h^{(i)}_{j} = 1$ for $j \neq i$.
Then $H$ possesses a relative ambiguity function in $G$ that is bounded above by a polynomial of degree $\abs{I}$.
\end{proposition}
\begin{proof}
For each $i \in I$ let $X_i$ be a finite generating set of $G_i$ and let $X = \bigcup \limits_{i \in I} \iota_i(X_i)$ denote the union of their canonical embeddings in $G$.
Using Lemma~\ref{LemmaChangeOfGenratingSet} we may assume that $X_i$ is contained in $Y_i$ for each $i \in I$.
We proceed with each factor $G_i$ as in Proposition~\ref{PropositionBoundedPreimage}.
To do so, we consider the kernel $N_i \defeq \ker((\psi_{I \setminus \{i\}})_H)$ for each $i \in I$ and its isomorphic image $\psi_i(N_i)$ in $G_i$.
Note that $\psi_i(N_i)$ contains the independent loxodromic elements $g^{(i)}_{i}$ and $h^{(i)}_{i}$ from our assumptions.
We can therefore apply
Remark~\ref{RemarkExtensionOfTheRelativeAmbFct}
which tells us that there is a concatenation map
\[
\Phi^{(i)} \colon G_i \times G_i \rightarrow G_i,\ (h,k) \mapsto h x^{(i)}_{h,k} k,
\]
with $x^{(i)}_{h,k} \in \psi_i(N_i)$ for all $h,k \in G_i$, that induces a linearly bounded ambiguity function $f_i$ of $G_i$.
By taking the product of the $\Phi^{(i)}$, we obtain a map
\[
\widetilde{\Phi} \colon G \times G \rightarrow G,\ ((h_i)_{i \in I},(h_i)_{i \in I}) \mapsto (h_i x^{(i)}_{h_i,k_i} k_i)_{i \in I}
\]
that restricts to a concatenation map $\Phi \colon H \times H \rightarrow H$.
We claim that $\Phi$ induces a polynomially bounded relative ambiguity function of $H$ with respect to $X$.
To see this, we fix an element $(g_i)_{i \in I} \in H$ and two numbers $m,n \in \N$.
We have to estimate the number of pairs $((h_i)_{i \in I},(k_i)_{i \in I}) \in B_{H}^{X}(m) \times B_{H}^{X}(n)$ that satisfy
\begin{equation}\label{EquationTupleOfEquations}
(h_i x^{(i)}_{h_i,k_i} k_i)_{i \in I} = (g_i)_{i \in I}.
\end{equation}
In view of Remark~\ref{RemarkUnionOfGenSets} it follows that the elements $h_i$ and $k_i$ in~\eqref{EquationTupleOfEquations} lie in $B_{G_i}^{X_i}(m)$ respectively in $B_{G_i}^{X_i}(n)$.
By the definition of $f_i$ it therefore follows that we the number of pairs $(h_i,k_i) \in B_{G_i}^{X_i}(m) \times B_{G_i}^{X_i}(n)$ that can appear in the $i$th coordinate in~\eqref{EquationTupleOfEquations} is bounded above by $f_i$.
%
%
Thus we see that the number of all $((h_i)_{i \in I},(k_i)_{i \in I}) \in B_{H}^{X}(m) \times B_{H}^{X}(n)$ that satisfy~\eqref{EquationTupleOfEquations} is bounded above by $f \defeq \prod \limits_{i \in I} f_i$, which is clearly bounded above by a polynomial of degree $\abs{I}$.
\end{proof}

\begin{theorem}\label{TheoremProducts}
Let $I$ be a finite set of indices.
For each $i \in I$ let $G_i$ be a finitely generated acylindrically hyperbolic group and let $G \defeq \prod \limits_{i \in I} G_i$ denote their product.
Let $H$ be a subgroup of $G$.
Suppose that $\psi_i(H) \leq G_i$ contains a generalized loxodromic element of $G_i$ for every $i \in I$.
Then $H$ has a relative ambiguity function in $G$ that is bounded above by a polynomial of degree $\abs{I}$.
Further the relative growth rate of $H$ exists with respect to every finite generating set of $G$.
\end{theorem}
\begin{proof}
For each $i \in I$ let $X_i$ be a finite generating set of $G_i$ and let $X = \bigcup \limits_{i \in I} \iota_i(X_i)$ denote the union of the canonical embeddings of the $X_i$ in $G$.
For each $i \in I$ we fix a generalized loxodromic element $h_i \in \psi_i(H)$.
Recall that Theorem~\ref{TheoremAcylActionOnCayleyGraph} provides us with generating sets $Y_i$ for $G_i$ such that $\Gamma(G_i,Y_i)$ is hyperbolic, $G_i$ acts acylindrically on $\Gamma(G_i,Y_i)$, and $h_i$ is a loxodromic element for this action.
By Lemma~\ref{LemmaDecomposition} there is a subset $I_1 \subset I$
such that the kernel of the restriction $(\psi_{I_1})_H$ is finite and such that the isomorphic image $\psi_i(N_i) \subset G_i$ of the kernel $N_i \defeq \ker((\psi_{I_1 \setminus \{i\}})_H)$ contains a loxodromic element $g_i \in G_i$ for every $i \in I_1$.
Due to Lemma~\ref{LemmaSmallKernel} and our choice of the generating set $X$ we may restrict our attention to the case where $I = I_1$.
In this case Theorem~\ref{TheoremClassificationAcylindric} tells us that $\psi_i(N_i)$
is either virtually cyclic or it contains infinitely many independent loxodromic elements.
Suppose that $\psi_i(N_i)$ is virtually cyclic for some $i \in I$.
Since $\psi_i(N_i)$ contains a loxodromic element, it follows that it is quasiisometrically embedded in $\Gamma(G_i,X_i)$.
Note that this implies that $N_i = \iota_i(\psi_{i}(N_i))$ is quasiisometrically embedded in $\Gamma(G,X)$.
From this it clearly follows that $\beta^X_{N_i}$ is bounded above by a linear function.
A further application of Lemma~\ref{LemmaSmallKernel} therefore allows us to replace $H$ by $\psi_{I \setminus \{i\}}(H)$ and $I$ by $I \setminus \{i\}$.
By an inductive application of this replacement we can assume that $\psi_{i}(N_i) \subset G_i$ contains two independent loxodromic elements for every $i \in I$.
Now the claim follows since this case was already treated in Proposition~\ref{PropositionTheSimpleCase}.
\end{proof}

\section{Right-angled Artin groups}\label{Applications}

In this section we will apply the results from the previous sections to study the relative growth rate in right-angled Artin groups.


\begin{definition}\label{DefinitionRAAG}
Let $\Gamma = \Gamma(V,E)$ be a finite graph with vertex set $V$ and edge set $E$.
The right-angled Artin group corresponding to $\Gamma$ is given by the presentation
\[
A_{\Gamma} = \langle \ V\ \vert \ [u,v] = 1\ \text{if} \ \{u,v\} \in E \ \rangle.
\]
\end{definition}

It is well-known that $A_{\Gamma}$ acts properly and cocompactly on a proper $\CAT(0)$-space $X_{\Gamma}$, which is given as the universal cover of the Salvetti complex associated to $A_{\Gamma}$.
This action can be used to distinguish generalized loxodromic elements in $A_{\Gamma}$.
In~\cite{Sisto18} Sisto considered the more general case of an arbitrary group $G$ that acts properly on a proper $\CAT(0)$-space $X$.
He showed that every element $g \in G$ that acts as a rank-$1$ isometry of $X$ lies in a hyperbolically embedded virtually cyclic subgroup of $G$.
In this case~\cite[Theorem~1.4]{Osin16} tells us that $g$ is a generalized loxodromic element of $G$.
In the case where $A_{\Gamma}$ acts on $X_{\Gamma}$, Behrstock and Charney showed that an element $g \in A_{\Gamma}$ is a rank-$1$ isometry if and only if $g$ cannot be conjugated into a join subgroup of $A_{\Gamma}$, i.e.\ a subgroup of the form $A_{\Gamma'}$ where $\Gamma'$ is a join in $\Gamma$ (see~\cite[Lemma~5.1]{BehrstockCharney12}).
For our purposes we need a criterion to decide whether a subgroup $H \leq A_{\Gamma}$ contains a rank-$1$ isometry.
The following theorem provides us with such a criterion in the case where $H$ is finitely generated.
In view of the above discussion we have replaced rank-$1$ isometries by generalized loxodromic elements.

\begin{theorem}[{\cite[Theorem~5.2]{BehrstockCharney12}\label{TheoremCriterionForContainingLoxo}}]
Let $H$ be a non-cyclic, finitely generated subgroup of $A_{\Gamma}$.
If $H$ cannot be conjugated into a join subgroup of $A_{\Gamma}$, then $H$ contains a rank-$1$ isometry.
\end{theorem}


The next result is a very useful tool in the study of right-angled Artin groups.
It was proven by Charney and Paris in the larger context of arbitrary Artin groups.

\begin{theorem}[{\cite[Theorem~1.2]{CharneyParis14}\label{TheoremConvexityOfParabolics}}]
Let $\Gamma = \Gamma(V,E)$ be a finite graph and let $\Gamma' = \Gamma(V',E')$ be a full subgraph of $\Gamma$.
Then $A_{\Gamma'}$ is convex in $A_{\Gamma}$, i.e.\ every geodesic in $\Gamma(A_{\Gamma},V)$ with endpoints in $A_{\Gamma'}$ is contained in $\Gamma(A_{\Gamma'},V')$.
\end{theorem}

We will use the following direct consequence of Theorem~\ref{TheoremConvexityOfParabolics}.

\begin{corollary}\label{CorollaryReductionEmptyGenerators}
Let $\Gamma = \Gamma(V,E)$ be a finite graph and let $\Gamma' = \Gamma(V',E')$ be a full subgraph of $\Gamma$.
If $H$ is a subgroup of $A_{\Gamma'}$, then $B_{H}^{V}(n) = B_{H}^{V'}(n)$ for every $n \in \N$.
Further a function $f \colon \N \rightarrow \N$ is a relative ambiguity function of $H$ with respect to $V$ if and only if $f$ is a relative ambiguity function of $H$ with respect to $V'$.
\end{corollary}

Corollary~\ref{CorollaryReductionEmptyGenerators} will enable us to restrict our attention to subgroups $H \leq A_{\Gamma}$ with the property that there is no vertex $v$ of $\Gamma = \Gamma(V,E)$ such that $H$ is contained in $V \setminus \{v\}$.

\begin{definition}\label{DefinitionSupportOfGroup}
Let $A_{\Gamma}$ be the right-angled Artin group corresponding to a finite graph $\Gamma = \Gamma(V,E)$ and let $H$ be a subgroup of $A_{\Gamma}$.
We define the \emph{support} of $H$ to be the subset of $V$ given by
$\supp(H) \defeq \{ v \in V\ \vert\ H \nsubseteq \langle V \setminus \{ v \} \rangle \}.$
\end{definition}

If the support $\supp(H)$ of a subgroup $H \leq A_{\Gamma}$ coincides with the vertex set of $\Gamma$ it may still happen that some conjugate $gHg^{-1}$ of $H$ has a strictly smaller support.
In the following it will be useful to avoid this situation.
In order to do so, we introduce the following lemma.



\begin{lemma}\label{LemmaAmbiguityForConjugates}
Let $G$ be a finitely generated group with a finite generating set $X$ and let $H$ be a subgroup of $G$.
Suppose that $H$ has a subexponential relative ambiguity function in $G$.
Then every conjugate $H^g$ of $H$ has a subexponential relative ambiguity function in $G$.
\end{lemma}
\begin{proof}
Let $f$ be a subexponential relative ambiguity function of $H$ with respect to some finite generating set $X$ of $G$.
Let further $\Phi \colon H \times H \rightarrow H,\ (h,k) \mapsto h x_{h,k} k$ be a concatenation map for $f$.
By replacing $X$ with $X^{g}$ and $\Phi$ with
\[
\Phi^g \colon H^g \times H^g \rightarrow H^g,\ (h^g,k^g) \mapsto h^g x_{h,k}^g k^g = (h x_{h,k} k)^g,
\]
we see that $f$ is a subexponential relative ambiguity function of $H$ in $G$ with respect to $X^g$.
\end{proof}

\begin{lemma}\label{LemmaReductionCyclicGenerators}
Let $A_{\Gamma}$ be a right-angled Artin group that decomposes as a direct product $A_{\Gamma} = A_{\Gamma_1} \times A_{\Gamma_2}$.
Let $H$ be a subgroup of $A_{\Gamma}$ such that the projection $\psi_2(H)$ is a cyclic subgroup of $A_{\Gamma_2}$.
If $\psi_1(H)$ has a subexponential relative ambiguity function in $A(\Gamma_1)$, then $H$ has a subexponential relative ambiguity function in $A_{\Gamma}$.
\end{lemma}
\begin{proof}
Suppose that $\psi_1(H)$ has a subexponential relative ambiguity function in $A(\Gamma_1)$.
By Lemma~\ref{LemmaSmallKernel} it suffices to show that the relative growth function of $\ker(\psi_1)$ in $A_{\Gamma_1}$ is subexponential.
As $\ker(\psi_1)$ canonically embeds into $\psi_2(H)$, we see that $\ker(\psi_1)$ is either trivial or an infinite cyclic group.
Without loss of generality we may assume the latter case.
Recall that right-angled Artin groups are $\CAT(0)$-groups and that cyclic subgroups of $\CAT(0)$-groups are quasiisometrically embedded (see~\cite[Theorem~III.$\Gamma$.1.1]{BridsonHaefliger99}).
Thus $\ker(\psi_1)$ is a quasiisometrically embedded cyclic subgroup of $A_{\Gamma}$, which implies that $\ker(\psi_1)$ has a linearly bounded, and hence subexponential, relative growth function in $A_{\Gamma}$.
\end{proof}

It is easy to see that finitely generated subgroups of free groups and free abelian groups are quasiisometrically embedded.
As right-angled Artin groups can be thought of as interpolating between free and free abelian groups, one might think that this is also true for them.
In this case it would seem quite clear that the relative growth rates of finitely generated subgroups of right-angled Artin groups should exist.
However, unlike for the free and free abelian case, there are right-angled Artin groups that contain exponentially distorted subgroups.
Examples of such groups that are discussed in~\cite{FloresKahrobaeiKoberda19}.
Nevertheless it turns out that the relative growth rate exists in these cases.

\begin{theorem}\label{TheoremExistenceForRAAGs}
Let $H$ be a finitely generated subgroup of a right-angled Artin group $A_{\Gamma}$.
Then the relative growth rate
of $H$ exists with respect to
every finite generating set of $A_{\Gamma}$.
\end{theorem}
\begin{proof}
Let $A_{\Gamma} = \prod \limits_{i \in I} A_{\Gamma_i}$ be the canonical direct product decomposition of $A_{\Gamma}$, let $V$ denote the vertex set of $\Gamma$, and let $V_i$ denote the vertex set of $\Gamma_i$ for each $i \in I$.
Recall that Corollary~\ref{CorollaryExistSubexpRelAmbiguity} tells us that it is sufficient to show that $H$ has a subexponential relative ambiguity function in $G$.
By Lemma~\ref{LemmaAmbiguityForConjugates} we can therefore replace $H$ by an arbitrary conjugate $H^g$.
Together with Corollary~\ref{CorollaryReductionEmptyGenerators} and Lemma~\ref{LemmaReductionCyclicGenerators} this allows us to assume that the support of every conjugate $H^g$ coincides with $V$ and that the image $\psi_i(H^g) \leq A_{\Gamma_i}$ is non-cyclic (and hence non-trivial) for every $i \in I$.
Note that this particularly tells us that $A_{\Gamma_i}$ is non-cyclic and directly indecomposable for every $i \in I$.
In view of Theorem~\ref{TheoremCriterionForContainingLoxo} this shows that the groups $A_{\Gamma_i}$ are acylindrically hyperbolic and
by Theorem~\ref{TheoremProducts} it therefore remains to prove that $\psi_i(H)$ contains a generalized loxodromic element of $A_{\Gamma_i}$. 
To see this we note that each $\psi_i(H)$, being a homomorphic image of $H$, is finitely generated.
We may therefore apply Theorem~\ref{TheoremCriterionForContainingLoxo}, which says that $\psi_i(H)$ contains a generalized loxodromic element if $\psi_i(H)$ is not conjugated into a join subgroup of $A_{\Gamma_i}$.
Suppose the contrary, i.e.\ that there is an index $j \in I$ and an element $g \in A_{\Gamma_j}$ such that $\psi_{j}(H)^g \leq A_{\Gamma_J}$ for some join subgraph $\Gamma_J = \Gamma_{J_1} \ast \Gamma_{J_2}$ of $\Gamma_j$. 
Since the direct factors $A_{\Gamma_i}$ are directly indecomposable, it follows that $\Gamma_J$ is a proper subgraph of $\Gamma_{j}$.
As $g$ lies in $A_{\Gamma_{j}}$, the other direct factors $A_{\Gamma_i}$ commute with $g$.
This implies that $\supp(H^g)$ lies in $\bigcup \limits_{i \in I \setminus \{j\}} V_i \cup J$, which is a proper subset of $V(\Gamma)$.
Now the theorem follows since this contradicts our assumptions on $H$.
\end{proof}


\textbf{\\Acknowledgments.} This article grew out of my master's thesis, which wrote at the Heinrich Heine University
under the supervision of Professor Oleg Bogopolski.
I would like to thank him for introducing me to the subject of hyperbolic groups and supporting me during my bachelor and master studies.

\bibliographystyle{amsalpha}
\bibliography{Literaturverzeichnis}

\providecommand{\bysame}{\leavevmode\hbox to3em{\hrulefill}\thinspace}
\providecommand{\MR}{\relax\ifhmode\unskip\space\fi MR }
\providecommand{\MRhref}[2]{%
  \href{http://www.ams.org/mathscinet-getitem?mr=#1}{#2}
}
\providecommand{\href}[2]{#2}
\begin{thebibliography}{CRSZ20}

\bibitem[BC12]{BehrstockCharney12}
Jason Behrstock and Ruth Charney, \emph{Divergence and quasimorphisms of
  right-angled {A}rtin groups}, Math. Ann. \textbf{352} (2012), no.~2,
  339--356. \MR{2874959}

\bibitem[BF10]{BestvinaFeighn10}
Mladen Bestvina and Mark Feighn, \emph{A hyperbolic {${\rm
  Out}(F_n)$}-complex}, Groups Geom. Dyn. \textbf{4} (2010), no.~1, 31--58.
  \MR{2566300}

\bibitem[BH99]{BridsonHaefliger99}
Martin~R. Bridson and Andr\'{e} Haefliger, \emph{Metric spaces of non-positive
  curvature}, Grundlehren der Mathematischen Wissenschaften [Fundamental
  Principles of Mathematical Sciences], vol. 319, Springer-Verlag, Berlin,
  1999. \MR{1744486}

\bibitem[Bow08]{Bowditch08}
Brian~H. Bowditch, \emph{Tight geodesics in the curve complex}, Invent. Math.
  \textbf{171} (2008), no.~2, 281--300. \MR{2367021}

\bibitem[CDS18]{CoulonDalBoSambusetti18}
R\'{e}mi Coulon, Fran\c{c}oise Dal'Bo, and Andrea Sambusetti, \emph{Growth gap
  in hyperbolic groups and amenability}, Geom. Funct. Anal. \textbf{28} (2018),
  no.~5, 1260--1320. \MR{3856793}

\bibitem[Coh82]{Cohen82}
Joel~M. Cohen, \emph{Cogrowth and amenability of discrete groups}, J. Funct.
  Anal. \textbf{48} (1982), no.~3, 301--309. \MR{678175}

\bibitem[CP14]{CharneyParis14}
Ruth Charney and Luis Paris, \emph{Convexity of parabolic subgroups in {A}rtin
  groups}, Bull. Lond. Math. Soc. \textbf{46} (2014), no.~6, 1248--1255.
  \MR{3291260}

\bibitem[CRSZ20]{CordesRussellSprianoZalloum20}
Matthew Cordes, Jacob Russell, Davide Spriano, and Abdul Zalloum,
  \emph{Regularity of morse geodesics and growth of stable subgroups}, 2020.

\bibitem[DFW19]{DahmaniFuterWise19}
Fran\c{c}ois Dahmani, David Futer, and Daniel~T. Wise, \emph{Growth of
  quasiconvex subgroups}, Math. Proc. Cambridge Philos. Soc. \textbf{167}
  (2019), no.~3, 505--530. \MR{4015648}

\bibitem[FKK19]{FloresKahrobaeiKoberda19}
Ram\'{o}n Flores, Delaram Kahrobaei, and Thomas Koberda, \emph{Algorithmic
  problems in right-angled {A}rtin groups: complexity and applications}, J.
  Algebra \textbf{519} (2019), 111--129. \MR{3874519}

\bibitem[Gri80]{Grigorchuk80b}
R.~I. Grigorchuk, \emph{Symmetrical random walks on discrete groups},
  Multicomponent random systems, Adv. Probab. Related Topics, vol.~6, Dekker,
  New York, 1980, pp.~285--325. \MR{599539}

\bibitem[Mil68]{Milnor68}
J.~Milnor, \emph{A note on curvature and fundamental group}, J. Differential
  Geometry \textbf{2} (1968), 1--7. \MR{0232311}

\bibitem[MM99]{MasurMinsky99}
Howard~A. Masur and Yair~N. Minsky, \emph{Geometry of the complex of curves.
  {I}. {H}yperbolicity}, Invent. Math. \textbf{138} (1999), no.~1, 103--149.
  \MR{1714338}

\bibitem[Ols17]{Olshanskii17}
A.~Yu. Olshanskii, \emph{Subnormal subgroups in free groups, their growth and
  cogrowth}, Math. Proc. Cambridge Philos. Soc. \textbf{163} (2017), no.~3,
  499--531. \MR{3708520}

\bibitem[Osi15]{Osin15}
D.~Osin, \emph{On acylindrical hyperbolicity of groups with positive first
  {$\ell^2$}-{B}etti number}, Bull. Lond. Math. Soc. \textbf{47} (2015), no.~5,
  725--730. \MR{3403956}

\bibitem[Osi16]{Osin16}
\bysame, \emph{Acylindrically hyperbolic groups}, Trans. Amer. Math. Soc.
  \textbf{368} (2016), no.~2, 851--888. \MR{3430352}

\bibitem[Rip82]{Rips82}
E.~Rips, \emph{Subgroups of small cancellation groups}, Bull. London Math. Soc.
  \textbf{14} (1982), no.~1, 45--47. \MR{642423}

\bibitem[Sha98]{Sharp98}
Richard Sharp, \emph{Relative growth series in some hyperbolic groups}, Math.
  Ann. \textbf{312} (1998), no.~1, 125--132. \MR{1645953}

\bibitem[Sis18]{Sisto18}
Alessandro Sisto, \emph{Contracting elements and random walks}, J. Reine Angew.
  Math. \textbf{742} (2018), 79--114. \MR{3849623}

\end{thebibliography}

\end{document}